\documentclass[letterpaper,12pt]{article}
\usepackage{slashbox}
\usepackage{latexsym,array,delarray,amsthm,amssymb,epsfig,newicktree}

\theoremstyle{plain}
\newtheorem{thm}{Theorem}

\newtheorem{prop}[thm]{Proposition}

\newtheorem{conj}[thm]{Conjecture}

\theoremstyle{definition}
\newtheorem{df}[thm]{Definition}
\newtheorem{ex}[thm]{Example}

\newcommand{\zz}{\mathbb{Z}}

\newcommand{\pp}{\mathbb{P}}

\newcommand{\rr}{\mathbb{R}}
\newcommand{\cc}{\mathbb{C}}

\begin{document}

\title{\bf The Mathematics of Phylogenomics}
\author{Lior Pachter and Bernd Sturmfels \\ Department of Mathematics,
UC Berkeley \\
\tt [lpachter,bernd]@math.berkeley.edu }
\date{\today }

\maketitle

\noindent  {\small ``The lack of real contact between mathematics and
biology is either a tragedy, a scandal or a challenge, it is hard to
decide which.'' -- Gian-Carlo Rota, \cite[p. 2]{DiscreteThoughts}}

\section{Introduction}

The grand challenges in biology today are being shaped by powerful
high-throughput technologies that have revealed the genomes
of many organisms, global expression patterns of genes and
detailed information about variation within populations. We are
therefore able to ask, for the first time, fundamental questions
about the evolution of genomes, the structure of genes and their
regulation,
and the connections between genotypes and
phenotypes of individuals.
The answers to these questions are all
predicated on progress in a variety of computational, statistical,
and mathematical fields \cite{Karp}.

The rapid growth in the characterization of genomes has
led to the advancement of a new discipline called
{\em Phylogenomics}. This discipline, whose scope and potential was
first outlined in \cite{E}, results from
the combination of two major fields in the life sciences:
{\em Genomics}, i.e., the study of the function and structure of genes
and genomes; and {\em Molecular Phylogenetics}, i.e., the study of the
hierarchical evolutionary relationships among organisms and their
genomes. The objective of this article is to offer
mathematicians a first introduction
to this emerging field, and to discuss specific
problems and developments arising from phylogenomics.

The mathematical tools to be highlighted 
in this paper are statistics, probability, combinatorics and -- last but not
least -- algebraic geometry.
Emphasis is placed on the use of
    {\em Algebraic Statistics}, which is the study of statistical
    models for discrete data using algebraic methods.
    See \cite[\S 1]{PSbook} for details.     Several models
    which are relevant for phylogenomics are shown to be
    algebraic varieties in certain high-dimensional
     spaces of probability distributions. This interplay between
     statistics and algebraic geometry offers a conceptual
     framework for
     (understanding existing and developing new) combinatorial algorithms
for biological sequence analysis. It is our hope that this  will  
contribute to
some ``real contact''  between mathematics  and biology.

    This paper is organized as follows. In Section 2
    we begin by reviewing the organization and structure of
   genomes. This section is meant as a  brief tutorial, aimed at
readers who have a little or no background in molecular biology.
It offers definitions of the relevant biological terminology.

Section 3 describes a very simple example of a statistical model for
inferring information about the genetic code. The point of this
example is to explain the philosophy of algebraic statistics:
{\em model means algebraic variety}.

A more realistic model, which is widely used in computational biology,
is the {\em hidden Markov model}. In Section 4 we explain this model
and discuss its applications to the problem of identifying genes in a genome.
Another key problem is the alignment of biological sequences. Section 5
reviews the statistical models and combinatorial algorithms for
{\em sequence alignment}. We also discuss the relevance of
{\em parametric inference} \cite{PS2} in this context.

In Section 6 we present statistical models for the evolution
of biological sequences. These models are algebraic varieties
associated with phylogenetic trees, and they play a key role
in inferring the ancestral relationships among  organisms
and in identifying regions in genomes that are under selection.

Section 7 gives an introduction to the field of {\em Phylogenetic
Combinatorics}, which is concerned with the combinatorics
and geometry of finite metric spaces, and their
application to data analysis in the life sciences.
We shall discuss  the space of all trees \cite{BHV},
the neighbor-joining algorithm for projecting
metrics onto this space,
and several natural generalizations of these concepts.

In Section 8 we go back to the data. We explain how one
obtains and studies DNA sequences generated by genome
sequencing centers, and we illustrate the mathematical
models by estimating the
probability that the DNA sequence in Conjecture \ref{MOLconj}
occurred by chance in ten vertebrate genomes.

\section{The Genome}

Every living organism has a genome, made up of deoxyribonucleic acids
(DNA) arranged in a double helix \cite{WC}, which
encodes  (in a way to be made precise) the fundamental ingredients of  
life.
Organisms are divided into two major classes: {\em eukaryotes}
(organisms whose cells contain a nucleus) and
{\em prokaryotes} (for example bacteria).
In our discussion we focus on genomes of eukaryotes,
and, in particular, the human genome \cite{Lander,Venter}.

Eukaryotic  genomes are divided into {\em chromosomes}. The human  
genome has
two copies of
each chromosome. There are 23 pairs of chromosomes:  22 {\em autosomes}
(two copies each in both men and women) and two {\em sex chromosomes},
which are denoted X and Y.
Women have two X chromosomes,
while men have one X and one Y chromosome.
Parents pass on a mosaic of their pair of chromosomes to their children.

The sequence of DNA molecules in a genome is typically represented as a
sequence of letters,
partitioned into chromosomes, from the four letter alphabet $\, \Omega = \{A,C,G,T\}$.
These letters correspond to the bases in the double helix, that is, the
{\em nucleotides} Adenine, Cytosine, Guanine and Thymine.
Since every base is paired with an
opposite base
($A$ with $T$, and $C$ with $G$ in the other half of the double helix),
in order to describe a genome it suffices to list the bases in only
one strand. However, it is important to note that the two strands have
a directionality which is indicated by the numbers
$5'$ and $3'$ on the ends (corresponding to carbon atoms in the helix
backbone). The convention is to represent DNA in the $5' \rightarrow
3'$ direction. The human genome
consists of approximately $2.8$ billion bases, and has been obtained
using high throughput sequencing
technologies that can be used to read the sequence of short DNA
fragments hundreds of bases long. {\em Sequence assembly algorithms}
are then used to piece together these fragments \cite{Myers}.
See also \cite[\S 4]{PSbook}.

\begin{table}
\small
\centering
\begin{tabular}{|c|c|c|c|c|} \hline
& \begin{tabular}{c}
T\\
\end{tabular}
  & \begin{tabular}{c}
C\\
\end{tabular}
  & \begin{tabular}{c}
A\\
\end{tabular}
  & \begin{tabular}{c}
G\\
\end{tabular}\\
\hline
\begin{tabular}{c}
T\\
\end{tabular}
& \begin{tabular}{l}
TTT $\mapsto$ Phe\\
TTC $\mapsto$ Phe\\
TTA $\mapsto$ Leu\\
TTG $\mapsto$ Leu\\
\end{tabular}
& \begin{tabular}{l}
TCT  $\mapsto$ Ser\\
TCC $\mapsto$ Ser\\
TCA $\mapsto$ Ser\\
TCG $\mapsto$ Ser\\
\end{tabular}
& \begin{tabular}{l}
TAT $\mapsto$ Tyr\\
TAC $\mapsto$ Tyr\\
TAA $\mapsto$ {\it stop} \\
TAG $\mapsto$ {\it stop} \\
\end{tabular}
& \begin{tabular}{l}
TGT $\mapsto$ Cys\\
TGC $\mapsto$ Cys\\
TGA $\mapsto$ {\it stop} \\
TGG $\mapsto$ Trp \\
\end{tabular}\\
\hline
\begin{tabular}{c}
C\\
\end{tabular}
& \begin{tabular}{l}
CTT $\mapsto$ Leu\\
CTC $\mapsto$ Leu\\
CTA $\mapsto$ Leu\\
CTG $\mapsto$ Leu\\
\end{tabular}
& \begin{tabular}{l}
CCT $\mapsto$ Pro\\
CCC $\mapsto$ Pro\\
CCA $\mapsto$ Pro\\
CCG $\mapsto$ Pro\\
\end{tabular}
& \begin{tabular}{l}
CAT $\mapsto$ His\\
CAC $\mapsto$ His\\
CAA $\mapsto$ Gln\\
CAG $\mapsto$ Gln\\
\end{tabular}
& \begin{tabular}{l}
CGT $\mapsto$ Arg\\
CGC $\mapsto$ Arg\\
CGA $\mapsto$ Arg\\
CGG $\mapsto$ Arg\\
\end{tabular}\\
\hline
\begin{tabular}{c}
A\\
\end{tabular}
& \begin{tabular}{l}
ATT $\mapsto$ Ile\\
ATC $\mapsto$ Ile\\
ATA $\mapsto$ Ile\\
ATG $\mapsto$ Met\\
\end{tabular}
& \begin{tabular}{l}
ACT $\mapsto$ Thr\\
ACC $\mapsto$ Thr\\
ACA $\mapsto$ Thr\\
ACG $\mapsto$ Thr\\
\end{tabular}
& \begin{tabular}{l}
AAT $\mapsto$ Asn\\
AAC $\mapsto$ Asn\\
AAA $\mapsto$ Lys\\
AAG $\mapsto$ Lys\\
\end{tabular}
& \begin{tabular}{l}
AGT $\mapsto$ Ser\\
AGC $\mapsto$ Ser\\
AGA $\mapsto$ Arg\\
AGG $\mapsto$ Arg\\
\end{tabular}\\
\hline
\begin{tabular}{c}
G\\
\end{tabular}
& \begin{tabular}{l}
GTT $\mapsto$ Val\\
GTC $\mapsto$ Val\\
GTA $\mapsto$ Val\\
GTG $\mapsto$ Val\\
\end{tabular}
& \begin{tabular}{l}
GCT $\mapsto$ Ala\\
GCC $\mapsto$ Ala\\
GCA $\mapsto$ Ala\\
GCG $\mapsto$ Ala\\
\end{tabular}
& \begin{tabular}{l}
GAT $\mapsto$ Asp\\
GAC $\mapsto$ Asp\\
GAA $\mapsto$ Glu\\
GAG $\mapsto$ Glu\\
\end{tabular}
& \begin{tabular}{l}
GGT $\mapsto$ Gly\\
GGC $\mapsto$ Gly\\
GGA $\mapsto$ Gly\\
GGG $\mapsto$ Gly\\
\end{tabular}\\
\hline
\end{tabular}
\caption{The genetic code.}
\end{table}

Despite the tendency
to abstract genomes as strings
over the alphabet $\Omega$, one must not forget that they
are highly structured: for example, certain subsequences within a
genome correspond to {\em genes}. These subsequences play the important
role of encoding {\em proteins}. Proteins are polymers made of twenty
different types of amino acids. Within a gene,
triplets of DNA, known as {\em codons}, encode the amino acids for the
proteins. This is
known as the {\em genetic code}. Table 1 shows the $64$ possible
codons, and the twenty amino acids they code for. Each amino acid is
represented by a three letter identifier
(``Phe'' = Phenylalanine, ``Leu'' = Leucin, ....).
   The three codons $TAA$, $TAG$ and $TGA$
    are special: instead of coding for an amino acid,
   they are used to indicate that the protein ends.

    In order to make protein, DNA is first copied into a similar molecule
called messenger RNA (abbreviated mRNA) in a process called {\em transcription}. It is the RNA that is {\em translated} into
protein. The entire process
    is referred to as {\em expression}.
    Proteins can be structural elements, or perform complex tasks (such  
as
regulation of expression) by interacting with the many
molecules and complexes in cells. Thus, the genome
is a blueprint for life. An understanding of the genes, the function of
their proteins, and their expression patterns is fundamental to
biology.

The human genome contains approximately $25,000$ genes, although the
exact number has still not been determined. While there
are experimental methods for validating and discovering genes, there is
still no known high throughput technology for accurately
identifying all the genes in a genome. The computational
problem of identifying genes, the
    {\em gene finding problem}, is an active area of research. One
of the main difficulties lies in the fact that only a small portion of
any genome is genic. For instance, less than $5 \%$ of the human genome
is known to be functional. In Section 4 we discuss this problem, and
the role of probabilistic models in formulating statistically sound
methods
for distinguishing genes from non-genic sequence. The models of choice,
hidden Markov models, allow for the integration of diverse biological
information (such as the genetic code and the structure of genes) and yet
are suitable for designing efficient algorithms.
By virtue of being algebraic
varieties, they provide a key example for the
link between algebra, statistics and genomics.
Nevertheless, the
current understanding of genes is not
sufficient to allow for the {\it ab-initio} identification of all the
genes in a genome, and it is through comparison with other genomes that
the genes are revealed \cite{ACP}.

The differences between the genomes of individuals in a population are
small and are primarily due to recombination events (part of the process by
which two copies of parental chromosomes are merged in the offspring).
On the other hand,
    the genomes of different {\em species} (classes of
    organisms that can produce offspring together)
tend to be much more divergent. Genome differences between species
can be explained by many biological events including:

\begin{itemize}
\item {\em Genome rearrangement} -- comparing  chromosomes of
related species reveals large segments that have been reversed and
flipped ({\em inversions}), segments that have been moved ({\em transpositions}), {\em fusions} of chromosomes, and other large scale
events. The underlying biological mechanisms are poorly understood  
\cite{PT, SN}.
\item {\em Duplications and loss} -- some genomes have undergone whole
genome duplications. This process was recently demonstrated for yeast
\cite{Kellis}. Individual chromosomes or genes may also be duplicated.
Duplication events are often accompanied by {\em gene loss}, as
redundant genes slowly lose or adapt their function over time  
\cite{EichS}.
\item {\em Parasitic expansion} -- large sections of genomes are
repetitive, consisting of
elements which can duplicate and re-integrate into a genome.
\item {\em Point mutation, insertion and deletion} -- DNA sequences
mutate, and in non-functional regions these mutations accumulate over time.
Such regions are also likely to exhibit deletions; for example,
strand slippage during replication can lead to an incorrect copy number
for repeated bases.
\end{itemize}
Accurate mathematical models for sequence
alignment and evolution, our topics in
Sections 5--7, have to take these processes into consideration.

\smallskip

Two distinct DNA bases that share a common ancestor are called {\em
homologous}. Homologous bases can be related via speciation and duplication 
events, and are therefore divided into two classes: {\em orthologous} and {\em paralogous}. Orthologous bases are descendant from a single base in an ancestral genome that underwent a speciation event, whereas two paralogous bases correspond to two distinct bases in a single ancestral genome that are related via a duplication. Because we cannot sequence ancestral
genomes, it is never possible to formally prove that two DNA bases are
homologous. However, statistical arguments can show that it is extremely likely
that two bases
are homologous, or even orthologous.
The problem of identifying homologous bases between genomes of related
species
is known as the {\em alignment problem}. We shall discuss this in
Section 5.

The alignment of genomes is the first step in identifying highly
conserved sequences that point to the small fraction of the genome
that is under selection, and therefore likely to be functional.
Although the problem of
sequence alignment is mathematically and computationally challenging,
proposed homologous sequences can be
rapidly and independently validated  (it is easy to check whether
two sequences align once they have been identified),
and the regions can often be tested
in a molecular biology laboratory to determine their function.
In other words, sequence alignment reveals concrete verifiable evidence
for evolutionary selection and often results in testable hypotheses.

As a focal point for our discussion, we present a specific DNA sequence 
of length $42$. This sequence was found in April 2004 as a byproduct of computational work 
conducted by Lior Pachter's group at Berkeley \cite{BP}.
Whole genome alignments were found and analyzed   of
the human (hs), chimpanzee (pt), mouse (mm), rat (rn), dog (cf),  
chicken (gg), frog (xt), zebra-fish (dr), fugu-fish (tr) and tetraodon  
(tn) genomes.
The abbreviations refer to the Latin names of these organisms. They will
be used in Table 3 and Figure \ref{NJtree}.
{}From alignments of the ten genomes, the following hypothesis
was derived, which we state in the form of a mathematical conjecture.

\begin{conj} \label{MOLconj}
    {\rm (The ``Meaning of Life'') \ }
The sequence of 42 bases
\begin{equation}
\label{MOL}
{\tt TTTAATTGAAAGAAGTTAATTGAATGAAAATGATCAACTAAG}
\end{equation}
was present in the genome of the ancestor of all vertebrates,
and it has been completely conserved to the present time
(i.e., none of the bases have been mutated,
nor have there been any insertions or deletions).
\end{conj}

The identification of such a sequence requires a
highly non-trivial computation: the alignment of ten genomes
(including mammalian genomes close to $3$ billion bases in length) and
subsequent analysis to identify conserved orthologous regions within the
alignment \cite{YP}. Using the tools described in Section 8, one checks
that the sequence (\ref{MOL}) is present in all ten genomes.
For instance, in the  human genome (May 2004 version), the sequence occurs on
chromosome 7 in positions 156501197--156501238. By
examining the alignment, one verifies that, with very high probability,
the regions containing this sequence in all ten genomes are orthologous.
Furthermore, the implied claim that (\ref{MOL}) occurs in all
present-day vertebrates can, in principle, be tested.

Identifying and analyzing sequences such as (\ref{MOL}) is important
because they are highly conserved yet often non-genic \cite{BPMSKMH}.
One of the ongoing mysteries in biology is to unravel the
function of the parts of the genome that is non-genic and
yet very conserved. The extent of
conservation points to the possibility of critical functions within the
genome. It may be a coincidence that the segment above contains two
copies of the {\em motif} \ {\tt TTAATTGAA}, but this motif may also
have some
function (for example, it may be bound by a protein). Indeed, the
identification of such elements is the first step towards
understanding the complex regulatory code of the genome. 
Back in 2003, we were amused to find that $42$ was the length of the longest
such sequence. In light of \cite{A},
it was decided to name this DNA-sequence ``The Meaning of Life''.

The conjecture was formulated in the spring of 2004
and it was circulated in the first {\tt arXiv} version of this paper.
In the fall of 2004, Drton, Eriksson and Leung \cite{DEL} conducted
a new study based on improved alignments. Their work 
led to even longer sequences with similar properties.
So, the Meaning of Life sequence no longer holds the
record in terms of length. However, since Conjecture \ref{MOLconj} 
has been inspiration for our group, and it still remains open today,
we decided to stick with this example.
It needs to be emphasized that disproving Conjecture \ref{MOLconj} 
would not invalidate any of the methodology presented
in this article. For a biological perspective we refer to \cite{DEL}.

\section{Codons}

Because of the genetic code, the set $\Omega^3$ of all three-letter
words over the alphabet $\, \Omega = \{A,C,G,T\}\,$
plays a special role in molecular biology.
As was discussed in Section 2, these words are called codons, with each
triplet
coding for one of $20$ amino acids (Table 1). The map from $64$ codons
to $20$ amino acids is not injective, and so multiple codons code for
the same amino acid. Such codons
are called {\em synonymous}. Eight amino acids have the property that the synonymous codons that code for them all agree in the first two positions. The third positions of such codons are called {\em four-fold degenerate}. 
The translation of a series of codons in a gene (typically a few
hundred)
results in a three-dimensional folded protein.

A {\em model for codons} is a statistical model
whose state space is the $64$-element set $\Omega^3$.
Selecting a model means specifying a family of
probability distributions $p = (p_{IJK})$
on $\Omega^3$. Each probability distribution $p$ is
    a $4 \times 4 \times 4$-table of non-negative
real numbers which sum to one. Geometrically,
a distribution on codons is a point $p$
in the $63$-dimensional probability simplex
$$  \Delta_{63} \quad = \quad \bigl\{ \,
\, p \in \rr^{\Omega^3} \,:\,
\sum_{IJK \in \Omega^3} p_{IJK} = 1 \,\,\,\,
\hbox{and} \,\,\,\, p_{IJK} \geq 0 \,\,\,
\hbox{for all} \, \,\, IJK \in \Omega^3\,\bigr\} . $$
A model for codons is hence nothing but a subset $\mathcal{M}$
of the simplex $\Delta_{63}$. Statistically meaningful
models are usually given in parametric form. If the number
of parameters is $d$, then there is a set $\mathcal{P} \subset \rr^{d}$
of allowed parameters, and the model $\mathcal{M}$ is the
image of a map $\phi$ from $\mathcal{P}$ into $\Delta_{63}$.
We illustrate this statistical point of view by means of a
very simple independence model.

Models for codons have played a prominent
role in the work of Samuel Karlin, who was one
of the mathematical pioneers in this field.
One instance of this is the {\em genome signature} in \cite{CMK}.
We refer to \cite[Example 4.3]{PSbook} for
a discussion of this model and 
more recent work on codon usage in genomes.

Consider a DNA sequence of length $3m$ which has been grouped
into $m$ consecutive codons. Let $u_{IJK}$ denote the number of
occurrences of a particular codon $IJK$. Then our data
is the $4 \times 4 \times 4$-table $u = (u_{IJK})$. The
entries of this table are non-negative integers, and
if we divide each entry by $m$ then we get
a new table $\frac{1}{m} \cdot u$ which is a point in the probability
    simplex $\Delta_{63}$.
This table is the {\em empirical distribution of codons}
in the given sequence.

Let $\mathcal{M}$ be the statistical model which stipulates that,
for the sequence under consideration, the first two positions in a codon
are independent from the third position. We may wish to test
whether this independence model fits our data $u$.
This question makes sense in molecular biology
    because many of the amino acids are
uniquely specified by the first two positions in any
codon which represents that particular amino acid (see Table 1).
Therefore, third
positions in synonymous codons tend to be independent of the first two.

Our independence model $\mathcal{M}$ has
$18$ free parameters. The set of allowed parameters
is an $18$-dimensional convex polytope, namely,
it is the product
$$ \mathcal{P} \quad = \quad \Delta_{15} \, \times \Delta_{3}.  $$
Here $\Delta_{15}$ is the $15$-dimensional simplex
consisting of probability distributions $\alpha = (\alpha_{IJ})$ on
$\Omega^2$, and $\Delta_3$ is the tetrahedron consisting
of probability distributions $\beta = (\beta_K) $ on $\Omega$.
Our model $\mathcal{M}$ is parameterized by the map
$$ \phi \,\,: \,\, \mathcal{P} \,\rightarrow \, \Delta_{63}\, , \quad
\phi((\alpha,\beta))_{IJK} \,= \,
\alpha_{IJ} \cdot \beta_K. $$
Hence $\,\mathcal{M} = {\rm image}(\phi) \,$ is an $18$-dimensional
algebraic subset inside  the $63$-dimensional simplex.
To test whether a given $4 \times 4 \times 4$-table $p$
lies in $\mathcal{M}$, we write that table as a two-dimensional
matrix with $16$ rows and $4$ columns:
$$
p' \quad = \quad
\pmatrix{
p_{AAA} & p_{AAC} & p_{AAG} & p_{AAT} \cr
p_{ACA} & p_{ACC} & p_{ACG} & p_{ACT} \cr
p_{AGA} & p_{AGC} & p_{AGG} & p_{AGT} \cr
p_{ATA} & p_{ATC} & p_{ATG} & p_{ATT} \cr
p_{CAA} & p_{CAC} & p_{CAG} & p_{CAT} \cr
    \vdots & \vdots & \vdots & \vdots \cr
    p_{TTA} & p_{TTC} & p_{TTG} & p_{TTT} \cr
}.
$$
Linear algebra  furnishes the following characterizations of
our model:

\begin{prop} For a point $p \in \Delta_{63}$,
the following conditions are equivalent:
\begin{enumerate}
\item The distribution $p$ lies in the model $\mathcal{M}$.
\item The $16 \times 4$ matrix $p'$ has rank $one$.
\item All $2 \times 2$-minors of the matrix $p'$ are zero.
\item
$ p_{IJK} \cdot p_{LMN} \,\, = \,\, p_{IJN} \cdot p_{L M K} \qquad
\hbox{for all nucleotides} \,\, I,J,K,L,M,N $.
\end{enumerate}
\end{prop}

In the language of  algebraic geometry, the
model $\mathcal{M}$ is known as the
{\em Segre variety}. More precisely, $\mathcal{M}$ is the
set of non-negative real points on the
Segre embedding of $\pp^{15} \times \pp^3$
in $\pp^{63}$. Here and throughout, the symbol
$\pp^{m}$ denotes the complex projective space
of dimension $m$.
One of the points argued in this
paper is that many of the more advanced
statistical  models, such as
graphical models \cite[\S 1.5]{PSbook}, actually used
in practice by computational
biologists are also algebraic varieties
with a special combinatorial structure.

Returning to our original biological motivation,
we are faced with the following statistics problem.
The DNA sequence under consideration is summarized
in the data $u$, and we wish to test whether or not
the model $\mathcal{M}$ fits the data. The geometric
idea of such a test is to determine whether or not the
empirical distribution $\,\frac{1}{m} \cdot u\,$
lies close to the Segre variety $\,\mathcal{M}$.
Statisticians have devised a wide range of such tests,
each representing a statistically
meaningful notion of  ``proximity to $\mathcal{M}$''.
These include the $\chi^2$-test, the $G^2$-test,
Fisher's exact test, and others, as explained
in standard statistics texts such as
\cite{Bic} or \cite{Fie}. A useful
tool of numerical linear algebra for
measuring the distance of a point
to the Segre variety is the
{\em singular value decomposition}
of the matrix $p'$. Indeed, $p'$ lies
on $\mathcal{M}$ if and only if
the second singular value of $p'$
is zero. Singular values
provide a good notion of distance between a given
    matrix and various determinantal varieties
such as $\mathcal{M}$.

One key ingredient in statistical tests is
    {\em maximum likelihood estimation}. The basic
idea is to find those model parameters
$\alpha_{IJ}$ and $\beta_K$ which would
best explain the observed data.
If we consider all possible genome sequences of length $3m$,
then the likelihood of observing our particular data $u$
equals
$$ \gamma \cdot \prod_{IJK \in \Omega^3}
p_{IJK}^{u_{IJK}}, $$
where $\gamma$ is a combinatorial constant.
This expression is a function of $(\alpha,\beta)$,
called the {\em likelihood function}. We wish to
find the point in our parameter domain
$\,\mathcal{P} \,=\, \Delta_{15} \times \Delta_3\,$
which maximizes this function. The solution
$(\hat \alpha, \hat \beta)$ to this non-linear
optimization problem is said to be the
{\em maximum likelihood estimate} for the data $u$.
In our independence model, the likelihood function is convex, and it is
easy
to write down the global maximum explicitly:
$$ \hat \alpha_{IJ} \,\, = \,\,
\frac{1}{m} \sum_{K \in \Omega} u_{IJK} \qquad
\hbox{and} \qquad
\hat \beta_K \,\, = \,\,
\frac{1}{m} \sum_{IJ \in \Omega^2} u_{IJK}. $$
In general, the likelihood function of a statistical
model will not be convex, and there is no easy
formula for writing the maximum likelihood estimate
as a function of the data. In practice, numerical
hill-climbing methods are used to
solve this optimization problem, but, of course,
there is no guarantee that a local maximum
found by such methods is actually the global maximum.

\section{Gene Finding}

In order to find genes in DNA sequences, it is necessary to identify
structural features and sequence
characteristics that distinguish genic sequence from non-genic
sequence. We begin by describing more of the detail of gene structure
which is essential in developing probabilistic models.

Genes are not contiguous subsequences of the genome,
but rather split into pieces called {\em introns} and {\em exons}.
After transcription, introns are spliced out and only the remaining
exons
are used in translation (Figure 1). Not all of the sequence
in the exons is translated; the initial and terminal
exons may consist of {\em untranslated regions} (indicated in grey in
the figure). Since the genetic code is in (non-overlapping) triplets,
it follows that the lengths of the translated portions of the exons
must sum to $0 \bmod 3$.
\begin{figure}[ht]
     \begin{center}
   \includegraphics[scale=0.7]{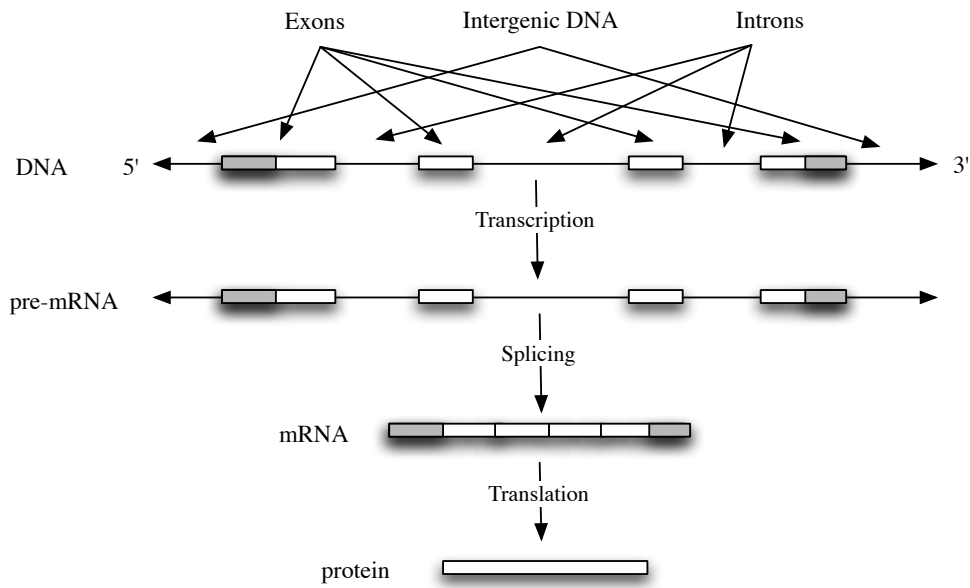}
     \end{center}
     \caption{Structure of a gene.}
     \label{fig:genestructure}
\end{figure}
In addition to the exon-intron structure of genes, there are
known sequence signals. The codon $ATG$ initiates translation, and thus
is the first codon following the untranslated portion of the initial
exons. The
final codon in a gene must be one of $TAG,TAA$ or $TGA$,
as indicted in Table 1.
    These codons  signal the translation machinery to stop. There are  
also
sequence signals at the intron-exon boundaries: $GT$ at the $5'$ end of
an intron and $AG$ at the $3'$ end.

A {\em hidden Markov model} (HMM) is a probabilistic model that allows
for simultaneous modeling of the bases in  a DNA sequence of length $n$
and the structural features associated with that sequence.
The  HMM consists of $n$ {\em observed}
random variables $Y_1,\ldots,Y_n$ taking on $l$ possible states,
and $n$ {\em hidden} random variables $X_1,\ldots,X_n$ taking on
$k$ possible states. In the context of phylogenomics, the
observed variables $Y_i$ usually have $l=4$ states, namely
$\Omega = \{A,C,G,T\}$. The hidden random variables  $X_i$
serve to model features associated with the sequence
which is generated by $Y_1, Y_2, \ldots, Y_n$. A
simple scenario is $k = 2$, with the
set of hidden states being
$\,\Theta \, = \, \{ \,{exon}, \,{intron} \}$.

The characteristic property of an HMM is that
the distributions of the $Y_i$ depend on the $X_i$,
while the $X_i$ form a {\em Markov chain}.
This is illustrated for $n=3$ in Figure 2, where
the unshaded circles represent the hidden  variables $X_1,X_2,X_3$
and the shaded circles represent the observed variables $Y_1,Y_2,Y_3$.

\begin{figure}[ht]
\begin{center}
 \includegraphics[scale=0.7]{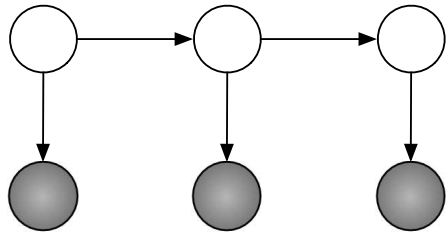}
     \end{center}
\caption{The hidden Markov model of length three.}
        \label{fig:HMM}\end{figure}

Computational biologists use HMMs to annotate DNA sequences.
The basic idea is this:  it is postulated that
the bases are instances of the random variables $Y_1,\ldots,Y_n$, and
the problem is to identify the most  likely assignments of states to
$X_1,\ldots,X_n$ that could be associated with the observations.
In gene finding, {\em homogeneous HMMs} are used.
This means that all transition probabilities
$\, X_i \rightarrow X_{i+1}\,$ are given
     by the same $k \times k$-matrix $S = (s_{ij})$,
and all the transitions $\,X_i \rightarrow Y_i \,$ are given by another
$k \times 4$-matrix $T = (t_{ij})$.
Here $s_{ij}$ represents the probability of
transitioning from hidden state $i$ to
hidden state $j$; for instance, if $k=2$ then
    $i,j \in \,\Theta \, = \, \{ \,{exon}, \,{intron} \}$.
The parameter $t_{ij}$ represents the probability
that state $i \in \Theta$ outputs  letter $j \in \Omega$.

In practice, the parameters $s_{ij}$ and $t_{ij}$
range over real numbers satisfying
\begin{equation}
\label{ParameterPolytope}
s_{ij}\,,\,\, t_{ij} \,\, \geq \,\, 0  \, \quad \hbox{and}\quad
    \sum_{j \in \Theta} s_{1j} \,\,\, =  \,\,\,
    \sum_{j \in \Omega} t_{1j} \,\,\, =  \,\,\, 1 .
\end{equation}
However, just like in our discussion of
the Segre variety in Section 3, we may relax the
requirements (\ref{ParameterPolytope})
and allow the parameters to be arbitrary complex
    numbers. This leads to the following
algebraic representation  \cite[\S 2]{PS1}.

\begin{prop} \label{HMMprop}
The homogeneous HMM is the image of a  map
    $\phi :  \cc^{k(k+l)} \rightarrow \cc^{l^n} \! $, where
each coordinate of $\phi$ is a
bi-homogeneous polynomial of degree  $n - 1$ in
the transition probabilities $s_{ij}$
and degree $n$ in the output probabilities $t_{ij}$.
\end{prop}

The coordinate $\phi_\sigma$ of the map $\phi$
indexed by a particular DNA sequence $\sigma \in \Omega^n$
represents the probability that the HMM
generates the  sequence $\sigma$.
The following explicit formula for that probability
establishes Proposition \ref{HMMprop}:
\begin{equation}
\label{BaumWelch}
\phi_\sigma \,\,\, = \,\,\,
\sum_{i_1 \in \Theta} t_{i_1 \sigma_1}
\biggl( \sum_{i_2 \in \Theta} s_{i_1 i_2} t_{i_2 \sigma_2}
\biggl( \sum_{i_3 \in \Theta} s_{i_2 i_3} t_{i_3 \sigma_3}
\biggl( \sum_{i_4 \in \Theta} s_{i_3 i_4} t_{i_4 \sigma_4}
\, \bigl( \,\, \cdots \,\, \bigr)
\biggr)
\biggr)
\biggr)
\end{equation}
The expansion of this polynomial has $k^n$ terms
\begin{equation}
\label{justoneterm}
t_{i_1 \sigma_1} s_{i_1 i_2} t_{i_2 \sigma_2}
    s_{i_2 i_3} t_{i_3 \sigma_3}  \cdots
    s_{i_{n-1} i_n} t_{i_n \sigma_n} .
\end{equation}
For any fixed parameters as in (\ref{ParameterPolytope}),
one wishes to determine a string $\,\widehat {\bf i}
    = (i_1,i_2,\ldots,i_n)\,\in\,\Theta^n \,$
which indexes a term  (\ref{justoneterm})
of largest numerical value among all $k^n$ terms of $\phi_\sigma$.
(If there is more than one string with maximum value
then we break ties lexicographically).
We call $\,\widehat {\bf i}\,$ the {\em explanation}
of the observation $\sigma$. In our example
$(k=2,l=4)$, the explanation  $\,\widehat {\bf i} \,$
of a DNA sequence $\sigma$ is
an element of $\,\Theta^n \, = \, \{ \,{exon}, \,{intron} \}^n$.
It reveals the crucial information of Figure 1, namely,
the location of the exons and introns.
In summary, the  DNA sequence to be annotated by an HMM
corresponds to the observation $\sigma \in \Omega^n$,
and the explanation $\widehat {\bf i}$ is the gene prediction.
Thus {\em gene finding}  means  nothing but
computing the output $\, \widehat {\bf i} \,$
from the input $\sigma$.

In real-world applications, the integer $n$ may be quite large.
It is not uncommon to annotate DNA sequences
of length $n \geq 1,000,000$.
The size  $k^n$ of  the search space
for finding the explanation is enormous (exponential in $n$).
Fortunately, the recursive decomposition in (\ref{BaumWelch}),
reminiscent of {\em Horner's Rule}, allows
us to evaluate a multivariate
polynomial with exponentially many
terms in linear time (in $n$).
In other words, for given numerical parameters $s_{ij}$ and $t_{ij}$,
we can compute the probability $\phi_\sigma(s_{ij},t_{ij})$
quite efficiently.

Similarly, the explanation $\widehat {\bf i}$
of an observed DNA sequence $\sigma$ can be computed in linear time.
This is done using the {\em Viterbi algorithm}, which
evaluates
$$
\max_{i_1 \in \Theta} T_{i_1 \sigma_1}
+ \biggl( \max_{i_2 \in \Theta} S_{i_1 i_2} + T_{i_2 \sigma_2}
+ \biggl( \max_{i_3 \in \Theta} S_{i_2 i_3} + T_{i_3 \sigma_3}
+ \biggl( \max_{i_4 \in \Theta} S_{i_3 i_4} + T_{i_4 \sigma_4}
    + \, \bigl(\, \cdots\, \bigr)
\biggr) \!
\biggr) \!
\biggr)
$$
where
$S_{ij} = {\rm log}(s_{ij})$ and
$T_{ij} = {\rm log}(t_{ij})$.
This expression is a piecewise linear
convex function on $\,\rr^{k(k+l)}$,
known as the {\em tropicalization}
of the polynomial $\phi_\sigma$.
Indeed, evaluating this expression requires
exactly the same operations as evaluating
$\phi_\sigma$, with the only difference
that we are replacing ordinary arithmetic
by the {\em tropical semiring}. The tropical
semiring (also known as the {\em max-plus} algebra) consists of the real numbers $\rr$ together
with an extra element $\infty$, where the arithmetic operations
of addition and multiplication are redefined to be $max$ (or equivalently $min$) and $plus$ respectively.
The tropical semiring and its use in dynamic programming optimizations is explained in \cite[\S 2.1]{PSbook}. 

Every choice of parameters $(s_{ij},t_{ij})$ specifies a
{\em gene finding function}
$$ \Omega^n \rightarrow \Theta^n ,\quad
\sigma \,\mapsto \, \widehat {\bf i} $$
which takes a sequence $\sigma$
to its explanation $\widehat {\bf i}$.
The number of all functions from
$\Omega^n$ to $\Theta^n$ equals
 $ \, 2^{n \cdot 4^n}\,$ and hence grows
 double-exponentially in $n$. However,
 the vast majority of these functions are  
 not gene finding functions. 
 The following remarkable complexity result
 was proved by Elizalde \cite{Eli}:

\begin{thm}
\label{exppol}
The number of gene finding functions grows
    at most polynomially in
the sequence length $n$.
\end{thm}

As an illustration consider the $n=3$
example visualized in Figure 2.
There are $\, 8^{64} = 6.277 \cdot 10^{57}\,$ functions
$\, \{A,C,G,T\}^3 \rightarrow \{ {exon}, {intron} \}^3\,$
but only a tiny fraction of these are gene finding functions.
(It would be interesting to determine the exact number).
It is an open problem to give a combinatorial characterization
of gene finding functions, and to come up with accurate
lower and upper bounds for their number as $n$ grows.

For gene finding HMMs, it is always the case that
$l$ is small and fixed (usually, $l=4$), and $n$ is large.
However, the size of $k$ or structure of
the state space for the hidden variables $X_i$
tends to vary a lot. While the  $k=2$ used in our discussion of
gene finding functions was meant to be just an illustration,
a biologically meaningful gene finding model could work with
just three hidden states: one for introns,
one for exons, and a state for
intergenic sequence.    However, in order to
enforce the constraint that the sum of the lengths of the
exons is $0 \bmod 3$, a more complicated
hidden state space is necessary.
Solutions to this problem were given in \cite{BK,KHRE}.

We conclude this section with a brief discussion of the
important problem of {\em estimating parameters} for HMMs.
Indeed, so far nothing has been said how the
values of the parameters $s_{ij}$ and $t_{ij}$
are to be chosen when running the Viterbi algorithm.
Typically, this choice involves a combination
of biological and statistical considerations.
Let us concentrate on the latter aspect.

Recall that maximum likelihood estimation is concerned with
finding  parameters for a statistical model 
which best explain the observed data.
As was the case for the codon model (Section 3), the
maximum likelihood estimate is an algebraic function of the data.
 In contrast to what we did
at the end of Section 3, it is now prohibitive
to locate the global maximum
in the polytope (\ref{ParameterPolytope}).
The {\em expectation-maximization} (EM) algorithm is  a
general technique used by statisticians to
find local maxima of the likelihood function \cite[\S 1.3]{PSbook}.
For HMMs, this algorithm is also known as
the {\em Baum-Welch algorithm}. It
takes advantage of the recursive decomposition in
(\ref{BaumWelch}) and it is  fast (linear in $n$). The widely used book \cite{DEKM} provides a 
good introduction to the use of the Baum-Welch algorithm in training HMMs for biological sequence 
applications. The connection between 
the EM algorithm and the Baum-Welch algorithm is explained in detail in \cite{HMY}.
In order to understand the performance of EM
or to develop more global methods \cite{CP}, it would be desirable
to obtain upper and lower bounds on the algebraic
degree \cite{HKS} of the maximum likelihood estimate.

\section{Sequence Alignment}

Although tools such as the hidden Markov model are  important for
modeling and analyzing individual genome sequences,
    the essence of phylogenomics lies in the power of {\em sequence
comparison}. Because functional sequences tend to accumulate fewer
mutations over time, it is possible, by comparing genomes, to identify
and characterize such sequences much more effectively.

    In this section we examine  models for sequence evolution that allow
for insertions,
    deletions and mutations in the special case of two genomes. These are
known as
    pairwise sequence alignment models. The specific model to be  
discussed
here
    is the {\em pair hidden Markov model}.
    In the subsequent section we shall examine
    phylogenetic models for more than two DNA sequences.

We have already seen two instances of statistical models that
are represented by polynomials in the model parameters
(the codon model and the hidden Markov model). Models for pairwise
sequence alignment are also specified by polynomials, and are in fact
close relatives of hidden Markov models.
What distinguishes the sequence alignment problem is an extra
layer of complexity which arises from a combinatorial explosion in the
number of possible alignments between sequences.
Here we describe one of the simplest alignment models (for a pair of
sequences), with a view towards connections with tree models and
algebraic statistics.

Given two sequences $\,\sigma^1 =
    \sigma^1_1 \sigma^1_2 \cdots \sigma^1_n \,$ and
$\,\sigma^2 =   \sigma^2_1 \sigma^2_2 \cdots \sigma^2_m \,$
over the alphabet $ \Omega = \{A,C,G,T \}$,
    an \emph{alignment} is a string over the auxiliary alphabet
$\{M,I,D\}$ such that
$\#M+\#D= n$ and $\#M+\#I=m $.
Here $\#M, \#I, \#D$ denote the number of characters $M,I,D$
in the word respectively.  An alignment records the ``edit steps'' from
the sequence
$\sigma^1$ to the sequence $\sigma^2$, where edit operations consist of
changing characters,
preserving them, or inserting/deleting them. An $I$ in the alignment
string
corresponds to an insertion from the first sequence to the second, a $D$ is a deletion
from the first
sequence to the second, and an $M$ is either a character change, or lack thereof.
The set ${\cal A}_{n.m}$  of all alignments depends only on the
integers $n$ and $m$, and not on  $\sigma^1$ and $\sigma^2$.

\begin{prop}
The cardinality of the set ${\cal A}_{n.m}$ of all alignments can be
computed
as the coefficient of the monomial $x^m y^n$ in the generating function
$$ \frac{1}{1-x-y-xy} \,\,\,\, = \,\,\,\,
1 + x+y +  x^2+3 xy+y^2 + \cdots +
x^5+9x^4 y+25 x^3y^2+ \cdots
$$
\end{prop}

    These cardinalities $| {\cal A}_{n,m}|$ are known as
    \emph{Delannoy numbers} in combinatorics
\cite[\S 6.3]{Stanley:99}. For instance, there are
$ \, | {\cal A}_{2,3}| \,= \, 25 \,$
alignments of two sequences of length two and three.
They are listed in Table 2 below.

\begin{figure}
     \begin{center}
   \includegraphics[scale=0.7]{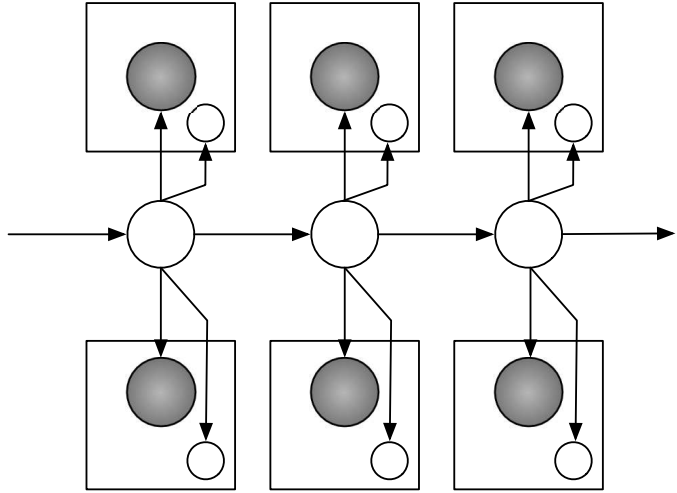}
     \end{center}
     \caption{A pair hidden Markov model for sequence alignment.}
     \label{fig:pairHMM}
\end{figure}

The {\em pair hidden Markov model} is visualized graphically
in Figure \ref{fig:pairHMM}.
The hidden random variables (unshaded nodes forming the Markov chain)
take on
the values $M,I,D$. Depending on the state at a hidden node,
either one or two characters are generated; in this way, pair hidden Markov models differ 
from standard hidden Markov models. The squares around the observed states (called plates)
are used to indicate that the number of characters generated may vary depending
on the hidden state. The number of characters generated is a random variable, 
indicated by unshaded nodes within the plates (called class nodes).
In pair hidden Markov models, the class nodes take on the values $0$ or $1$ corresponding to whether
or not a character is generated.
Pair hidden Markov models are therefore HMMs, where the structure of
the model depends on the assignments to the hidden states. The graphical
model structure of pair HMMs is explained in more detail in \cite{ABP}.

The next proposition gives the algebraic representation
of the pair hidden Markov model.
For a given alignment ${\bf a} \in {\cal A}_{n,m}$, we
denote the $j$th character in $\,{\bf a}\,$ by $a_j$, we
write $\,a[i] \,$ for $\,\#M+\#D \,$ in the prefix
$\,a_1 a_2 \ldots a_i$, and we write
$\,a \langle j \rangle\,$ for $\,\#M+\#I\,$ in
the prefix $\,a_1 a_2 \ldots a_j$.
Let $\sigma^1$ and $\sigma^2$ be two DNA sequences of lengths $n,m$
respectively.
Then the probability that our model generates these two sequences equals
\begin{equation}
\label{pairhmm}
\phi_{\sigma^1,\sigma^2}
    \quad = \quad \sum_{{\bf a} \in {\cal A}_{n,m}}
    t_{a_1}(\sigma^1_{a[1]},\sigma^2_{a \langle 1 \rangle}) \cdot
\prod_{i = 2}^{|{\bf a}|}
s_{a_{i-1} a_i} \cdot t_{a_i}(\sigma^1_{a[i]},\sigma^2_{a \langle i
\rangle}) ,
\end{equation}
where the parameter $\,s_{a_{i-1} a_i}\,$ is the transition probability
from state $a_{i-1}
$ to $a_i$, and the parameter $\,t_{a_i}(\sigma^1_{a[i]},\sigma^2_{a
\langle i \rangle})\,$
is the output probability
for a given state $a_i$ and the indicated output characters on the
strings
$\sigma^1$ and $\sigma^2$.

\begin{prop} \label{pairHMMmap}
The pair hidden Markov model for sequence alignment is the
image of a polynomial map $\,\phi :  \cc^{33}
\rightarrow \cc^{4^{n+m}}$.
The coordinates of the map $\phi$ are
the polynomials of degree $\leq 2n  + 2m - 1 $ which are
given in  (\ref{pairhmm}).
\end{prop}

We need to explain why the number of parameters in
our representation of the
pair hidden Markov model  is $33$.
First, there are nine parameters
$$ S \quad = \quad
\pmatrix{
s_{MM} &  s_{MI} &  s_{MD} \cr
s_{IM} &  s_{II} &  s_{ID} \cr
s_{DM} &  s_{DI} &  s_{DD}
} $$
which play the same role as in Section 4,
namely, they represent transition probabilities
in the Markov chain. There are
$16$ parameters $\,t_M(a,b) =: t_{Mab}\,$
for the probability that letter $a$ in
$\sigma^1$ is matched with letter $b$ in $\sigma^2$.
The insertion parameters $\,t_I(a,b) \,$
depend only on the letter $b$, and the
deletion parameters  $\,t_D(a,b) \,$
depend only on the letter $a$, so there
are only  $8 $ of these parameters.
Hence the total number of (complex) parameters is
$\,9 + 16 + 8 \,=\, 33$. Of course, in our applications,
probabilities are non-negative reals that sum to one,
so we get a reduction in the number of parameters,
just like in (\ref{ParameterPolytope}).
    In the upcoming example,
which explains the algebraic representation of
Proposition \ref{pairHMMmap},
we use the abbreviations $\,t_{Ib}\,$ and $\,t_{Da}\,$
for these parameters.

\begin{table}
\centering
\begin{tabular} {|l|l|l|}  \hline
IIIDD & \,\, $( \,\cdot \cdot \cdot ij \,,\, klm\cdot \cdot  \, )$ & $
    t_{Ik} s_{II} t_{Il} s_{II} t_{Im} s_{ID} t_{Di} s_{DD} t_{Dj} $\\
IIDID & \,\, $( \,\cdot \cdot i\cdot j \, ,\, kl\cdot m\cdot  \, )$ & $
    t_{Ik} s_{II} t_{Il} s_{ID} t_{Di} s_{DI} t_{Im} s_{ID} t_{Dj} $\\
IIDDI & \,\, $( \,\cdot \cdot ij \,\cdot \,,\, kl\cdot \cdot m \, )$ & $
    t_{Ik} s_{II} t_{Il} s_{ID} t_{Di} s_{DD} t_{Dj} s_{DI} t_{Im} $\\
IDIID & \,\, $( \,\cdot \, i\cdot \cdot j\,,\, k\cdot lm\cdot  \, )$ & $
    t_{Ik} s_{ID} t_{Di} s_{DI} t_{Il} s_{II} t_{Im} s_{ID} t_{Dj} $\\
IDIDI & \,\, $( \,\cdot \, i\cdot j\cdot \,,\, k\cdot l\cdot m \, )$ & $
    t_{Ik} s_{ID} t_{Di} s_{DI} t_{Il} s_{ID} t_{Dj} s_{DI} t_{Im} $\\
IDDII & \,\, $( \,\cdot \,ij \cdot \cdot \,,\, k\cdot \cdot lm \, )$ & $
    t_{Ik} s_{ID} t_{Di} s_{DD} t_{Dj} s_{DI} t_{Il} s_{II} t_{Im} $\\
DIIID & \,\, $( \,i\cdot \cdot \cdot j \,,\, \cdot \, klm\cdot  \, )$ &
$
    t_{Di} s_{DI} t_{Ik} s_{II}
t_{Il} s_{II} t_{Im} s_{ID} t_{Dj} $\\
DIIDI & \,\, $( \,i\cdot \cdot j\cdot \,,\, \cdot \,kl\cdot m \, )$ & $
    t_{Di} s_{DI} t_{Ik} s_{II} t_{Il} s_{ID} t_{Dj} s_{DI} t_{Im} $\\
DIDII & \,\, $( \,i\cdot j\cdot \cdot \,,\, \cdot \,k\cdot lm \, )$ & $
    t_{Di} s_{DI} t_{Ik} s_{ID} t_{Dj} s_{DI} t_{Il} s_{II} t_{Im} $\\
DDIII & \,\, $( \,ij\cdot \cdot \,\cdot\, ,\, \cdot \cdot klm \, )$ & $
    t_{Di} s_{DD} t_{Dj} s_{DI} t_{Ik} s_{II} t_{Il} s_{II} t_{Im} $\\
MIID & \,\, $( \,i\cdot \cdot j \,,\, klm \,\cdot  \, )$ & $
t_{Mik} s_{MI} t_{Il} s_{II} t_{Im} s_{ID} t_{Dj} $\\
MIDI & \,\, $( \,i\cdot j\cdot \,,\, kl\cdot m \, )$ & $     t_{Mik}
s_{MI} t_{Il} s_{ID} t_{Dj} s_{DI} t_{Im} $\\
MDII & \,\, $( \,ij\cdot \cdot \,,\, k\cdot lm \, )$ & $     t_{Mik}
s_{MD} t_{Dj} s_{DI} t_{Il} s_{II} t_{Im} $\\
IMID & \,\, $( \,\cdot \,i\cdot j \,,\, klm\cdot  \, )$ & $     t_{Ik}
s_{IM} t_{Mil} s_{MI} t_{Im} s_{ID} t_{Dj} $\\
IMDI & \,\, $( \,\cdot \,ij\,\cdot \,,\, kl\cdot m \, )$ & $     t_{Ik}
s_{IM} t_{Mil} s_{MD} t_{Dj} s_{DI} t_{Im} $\\
IIMD & \,\, $( \,\cdot \cdot ij\,,\, klm \,\cdot  \, )$ & $     t_{Ik}
s_{II} t_{Il} s_{IM} t_{Mim} s_{MD} t_{Dj} $\\
IIDM & \,\, $( \,\cdot \cdot ij\,,\, kl\cdot m \, )$ & $     t_{Ik}
s_{II} t_{Il} s_{ID} t_{Di} s_{DM} t_{Mjm} $\\
IDMI & \,\, $( \,\cdot ij\cdot \,,\, k\cdot lm \, )$ & $     t_{Ik}
s_{ID} t_{Di} s_{DM} t_{Mjl} s_{MI} t_{Im} $\\
IDIM & \,\, $( \,\cdot i\cdot j\,,\, k\cdot lm \, )$ & $     t_{Ik}
s_{ID} t_{Di} s_{DI} t_{Il} s_{IM} t_{Mjm} $\\
DMII & \,\, $( \,ij\cdot \cdot \,,\, \cdot \,klm \, )$ & $     t_{Di}
s_{DM} t_{Mjk} s_{MI} t_{Il} s_{II} t_{Im} $\\
DIMI & \,\, $( \,i\cdot j\cdot \,,\, \cdot \,klm \, )$ & $     t_{Di}
s_{DI} t_{Ik} s_{IM} t_{Mjl} s_{MI} t_{Im} $\\
DIIM & \,\, $( \,i\cdot \cdot j\,,\, \cdot \,klm \, )$ & $     t_{Di}
s_{DI} t_{Ik} s_{II} t_{Il} s_{IM} t_{Mjm} $\\
MMI & \,\, $( \,ij \,\cdot\,\, , \,\,klm \, )$ & $     t_{Mik} s_{MM}
t_{Mjl} s_{MI} t_{Im} $\\
MIM & \,\, $( \,i \cdot j \,\,,\,\, klm \, )$ & $     t_{Mik} s_{MI}
t_{Il} s_{IM} t_{Mjm} $\\
IMM & \,\, $( \,\cdot \,ij \,\,,\,\, klm \, )$ & $     t_{Ik} s_{IM}
t_{Mil} s_{MM} t_{Mjm} $\\ \hline
\end{tabular}
\caption{Alignments for a pair of sequences of length $2$ and $3$.}
\end{table}

Consider two sequences $\, \sigma^1 = ij \,$ and $\sigma^2 = klm \,$
of length $n = 2$ and $m = 3$ over the alphabet $\Omega = \{A,C,G,T \}$.
The number of alignments is $\,| {\cal A}_{2,3} | = 25$, and they are
listed in Table 2. For instance, the alignment
$\, MIID $, here written $( \,i\cdot \cdot j \,,\, klm \,\cdot  \, )$,
corresponds to
$\, \matrix{ {\tt i} \!\!\! & \!\!\! {\tt -} \!\!\! & \!\!\! {\tt -} 
\!\!\! & \!\!\! {\tt j} \cr
             {\tt k} \!\!\! & \!\!\! {\tt l} \!\!\! & \!\!\! {\tt m} 
\!\!\! & \!\!\! {\tt -}}\,$ 
in standard genomics notation.

The polynomial $\, \phi_{\sigma^1,\sigma^2}\, $ is the sum of the
$25$ monomials (of degree $9,7,5$) in the rightmost column.
Thus the pair hidden Markov model
 presented in Table 2 is nothing but a polynomial map
$$ \phi \,\, : \,\, \cc^{33} \,\,\rightarrow \,\, \cc^{1024}. $$

Statistics is all about making inferences. We shall now explain how
this is done with this model.
For any fixed parameters $\,s_{\cdot \cdot} \,$ and $\,t_{\cdot
\cdot}\,$,
one wishes to determine the alignment
$\,\widehat {\bf a} \in {\cal A}_{n,m}\,$ which indexes
the term of largest numerical value among the
Delannoy many terms of the polynomial $\phi_{\sigma_1,\sigma_2}$.
(If there is more than one alignment with maximum value then
we break ties lexicographically). We call  $\,\widehat {\bf a} \,$ the
{\em explanation} of the observation $\,(\sigma^1,\sigma^2)$.

The explanation for a pair of DNA sequences can be computed in
polynomial time (in their lengths $n$ and $m$) using a variant
of the Viterbi algorithm. Just like in the previous section,
the key idea is to {\em tropicalize}
the coordinate polynomials (\ref{pairhmm})
of the statistical model in question.  Namely, we compute
\begin{equation}
\label{supportfct}
\max_{{\bf a} \in {\cal A}_{n,m}} \,\,
    T_{a_1}(\sigma^1_{a[1]},\sigma^2_{a \langle 1 \rangle}) +
\sum_{i = 2}^{|{\bf a}|}
S_{a_{i-1} a_i} + T_{a_i}(\sigma^1_{a[i]},\sigma^2_{a \langle i
\rangle}) ,
\end{equation}
where $ \, S_{\cdot \cdot} \, = \,{\rm log}(s_{\cdot \cdot}) \,$
and $ \, T_{\cdot \cdot} \, = \,{\rm log}(t_{\cdot \cdot}) $.
The ``arg max'' of this piecewise linear convex function
is the optimal alignment $\,\widehat {\bf a}$.
Inference in the pair HMM means computing the
optimal alignment of two observed DNA sequences.
In other words, by inference we mean evaluating the {\em alignment
function}
$$ \Omega^n \times \Omega^m \,\, \rightarrow \,\, {\cal A}_{n,m} \, ,
\quad
(\sigma^1, \sigma^2) \, \mapsto \, \widehat {\bf a}. $$
There are doubly-exponentially many functions
from $\Omega^n \times \Omega^m$ to ${\cal A}_{n,m}$, but, by 
Elizalde's {\em Few Inference Functions Theorem} \cite{Eli},
 at most polynomially many of them are
alignment functions. 
Like for gene finding functions
(cf.~Theorem \ref{exppol}),
it is an open problem to characterize
alignment functions.

The function $\, \rr^{33} \rightarrow \rr \,$
given in (\ref{supportfct}) is the support function of a convex polytope
in $\rr^{33}$, namely, the {\em Newton polytope} of the polynomial
$\,\phi_{\sigma^1, \sigma^2}$. The vertices of this polytope correspond
to all optimal alignments of the sequences
$\sigma^1, \sigma^2$, with respect to  all possible
choices of the parameters, and the normal fan of the polytope divides
the logarithmic parameter space into regions which yield the same optimal
alignment. This can be used for analyzing the sensitivity of alignments to
parameters, and for the computation of posterior probabilities of
optimal alignments. The process of computing
this polytope is called {\em parametric alignment}
or {\em parametric inference}. It is known
\cite{FV, PS2, Waterman:92}
that parametric inference can be done in polynomial
time (in $m$ and $n$).

An important remark is that the formulation of sequence alignment with
pair Hidden Markov models is equivalent to
combinatorial ``scoring schemes'' or ``generalized edit distances''
which
can be used to assign weights to alignments \cite{BH}.
The simplest scoring scheme consists of two parameters:
    a mismatch score $mis$, and an indel score $gap$ \cite{Gusfield:97}.
The weight of an alignment is the sum of the scores for all positions
in the alignment, where a match gets a score of $1$.
In the case where $mis$ and $gap$ are non-negative, this is equivalent to specializing the $33$ logarithmic parameters
$ \, S_{\cdot \cdot} \, = \,{\rm log}(s_{\cdot \cdot}) \,$
and $ \, T_{\cdot \cdot} \, = \,{\rm log}(t_{\cdot \cdot}) \,$
of the pair hidden Markov model as follows:
\begin{eqnarray*}
& S_{ij} = 0,  \quad
T_{Ij} =  T_{Di} =- gap\  \hbox{for all $i,j$,} \\
& T_{Mij}=-1 \hbox{ if $i=j$}, \,\,\, \hbox{and} \,\,\,
    T_{Mij}=-mis\, \hbox{ if $i \neq j$}.
\end{eqnarray*}
The case where the scoring scheme consists of both positive and negative parameters corresponds
to a normalized pair hidden Markov model \cite{DEKM}.
This specialization of the parameters
corresponds to projecting the Newton polytope
of $\phi_{\sigma^1,\sigma^2}$ into two dimensions.
Parametric alignment means computing the
resulting two-dimensional polygon. For two sequences of length $n$, an
upper bound on the number of vertices in the polygon is
$O(n^{2/3})$.  We have observed that for  biological sequences
the number may be much smaller.
See \cite{FV} for a survey from the perspective of
computational geometry.

In the strict technical sense, our polynomial
formulation (\ref{pairhmm}) is not needed to
derive or analyze combinatorial algorithms
for sequence alignment. However,  the translation from
algebraic geometry (\ref{pairhmm})
to discrete optimization (\ref{supportfct})
offers much more than just esthetically pleasing formulas.
We posit that (tropical) algebraic geometry is
a conceptual framework for developing
new models and designing new algorithms
of practical value for phylogenomics.

\section{Models of Evolution}

Because organisms from
different species cannot produce offspring together, mutations and
genome changes that occur within a species
are independent of those occurring in another species.
There are some exceptions to this statement, such
as the known phenomenon of {\em horizontal transfer} 
in bacteria which results in the transfer of genetic material 
between different species, however we
ignore such scenarios in this discussion.
We can therefore represent the evolution of species (or phyla) via a
tree structure. The study of tree structures
in genome evolution is referred to as {\em phylogenetics}.
A phylogenetic {\em $X$-tree} is a tree $T$ with all internal vertices
of degree at least 3, and with the leaves labeled by a set $X$ which
consists of different species. In this section, we assume 
that $T$ is known and that vertices in $T$ correspond to 
known speciation events. We begin by describing statistical 
models of evolution that are used to identify regions
between genomes that are under selection.

Evolutionary models attempt to capture three important aspects of
evolving sequences: {\it branch length}, {\it substitution} and {\it
mutation}. Consider a single ancestral base $b$ at the root $r$ of a
phylogenetic tree $T$, and assume that there are no insertions or
deletions over time.
Since the ancestral base changes, it is possible that at two leaves
$x,y \in X$ we observe bases $c_1 \neq c_2$. We say that there has been
a {\em substitution} between $x$ and $y$. In a probabilistic model of
evolution, we would like to capture the possibility for change along
internal edges of the tree, with the possibility of back substitutions
as well. For example, it is possible that $b \rightarrow c_1
\rightarrow b \rightarrow c_1$ along the path from $r$ to $x$.

\begin{df}
A {\em rate matrix} (or {\em $Q$-matrix}) is a square matrix
     $Q=(q_{ij})_{i,j \in \Omega}$ (with rows and columns
     indexed by the nucleotides) satisfying the properties
\[ q_{ij} \geq 0 \quad \,\, \hbox{for} \,\,\, \ i \neq j, \]
\[ \sum_{j \in \Omega} q_{ij} =0 \quad \,\,\hbox{for all} \,\, \ i \in
\Omega, \]
\[ q_{ii}  < 0 \quad \hbox{for all} \,\, \ i \in \Omega. \]
\end{df}

Rate matrices capture the notion of {\em instantaneous rate of
mutation}. {}From a given rate matrix $Q$ one computes the {\em  
substitution
matrices} $P(t)$ by exponentiation.
The entry of $P(t)$ in row $b $ and column $c$ equals the  probability
that the substitution $\,b \rightarrow \cdots \rightarrow c\,$ occurs
in a time interval of length $t$.
We recall the following well-known result
about continuous-time Markov models.

\begin{prop}
Let $Q$ be any rate matrix and $\,P(t)\,=\,e^{Qt} =
\sum_{i=0}^\infty \frac{1}{i \, !} Q^i t^i  $. Then
\begin{enumerate}
\item $ P(s+t)\,=\,P(s)+P(t) $,
\item $ P(t)\,$ is the unique solution to $\, P'(t)=P(t) \cdot Q,
\,P(0)= {\bf 1} \,$ for $\,t \geq 0$,
\item $P(t)\,$ is the unique solution to $\, P'(t) = Q \cdot P(t),\,
P(0)= {\bf 1} \,$ for $\,t \geq 0$.
\end{enumerate}
Furthermore, a matrix $Q$ is a rate matrix if and only if
the matrix $P(t) = e^{Qt}$ is a stochastic matrix (nonnegative
with row sums equal to one) for every $t$.
\end{prop}

The simplest model is the {\em Jukes-Cantor DNA model}, whose rate matrix
is
\[ Q \quad = \quad
\pmatrix{
-3 \alpha & \alpha & \alpha & \alpha \cr
\alpha & -3 \alpha & \alpha & \alpha \cr
\alpha & \alpha & -3 \alpha & \alpha \cr
\alpha & \alpha & \alpha & -3 \alpha
}  , \]
where $\alpha \geq 0 $ is a parameter. The corresponding substitution
matrix equals
\[ P(t) \quad = \quad \frac{1}{4}
\pmatrix{
1+3e^{-4 \alpha t} & 1-e^{-4\alpha t} & 1-e^{-4\alpha t}  &
1-e^{-4\alpha t} \cr
1-e^{-4\alpha t} & 1+3e^{-4 \alpha t} & 1-e^{-4\alpha t}  &
1-e^{-4\alpha t} \cr
1-e^{-4 \alpha t} & 1-e^{-4\alpha t} & 1+3e^{-4 \alpha t} &
1-e^{-4\alpha t} \cr
1-e^{-4 \alpha t} & 1-e^{-4\alpha t} & 1-e^{-4\alpha t} & 1+3e^{-4
\alpha t}  \cr
}  . \]
The expected number of substitutions over time $t$ is the quantity
\begin{equation}
\label{branchlength} 3 \alpha t \quad = \quad -\frac{1}{4} \cdot {\rm
trace}(Q) \cdot t
     \quad = \quad - \frac{1}{4} \cdot {\rm log} \,{\rm det} \bigl( P(t)
\bigr).
     \end{equation}
     This number is called the {\em branch length}. It can be computed
from
     the substitution matrix $P(t)$ and is used
     to weight the edges in a phylogenetic $X$-tree.

     One way to specify an {\em evolutionary model} is to give a
phylogenetic $X$-tree $T$ together with a rate matrix $Q$
     and an initial distribution for the root of $T$ (which we
     here assume to be the stationary distribution on $\Omega$).
     The branch lengths of the edges are unknown parameters,
     and the objective  is to estimate these branch lengths from
      data. Thus if the tree $T$ has $r$ edges, then such a  model
     has $r$ free parameters, and, according to the philosophy
     of algebraic statistics, we would like to regard it as an
     $r$-dimensional algebraic variety.

Such an algebraic representation does indeed exist. We shall
explain it for the Jukes-Cantor DNA model on an
$X$-tree $T$. Suppose that $T$ has $r$ edges and $|X| = n$ leaves.
Let $P_i(t)$ denote the substitution matrix
associated with the $i$-th edge of the tree.
We write $\, 3 \alpha_i t_i \, = \,
- \frac{1}{4} {\rm log} \,{\rm det} \bigl( P_i(t) \bigr)\,$
for the branch length of the $i$-th edge, and we set
$\pi_i = \frac{1}{4} (1-e^{-4 \alpha_i t_i})$ and
$\theta_i = 1 - 3  \pi_i$. Thus
$$ P_i(t) \quad = \quad
\pmatrix{\theta_i & \pi_i & \pi_i & \pi_i  \cr
\pi_i & \theta_i & \pi_i & \pi_i   \cr
\pi_i & \pi_i &\theta_i & \pi_i   \cr
\pi_i & \pi_i & \pi_i &\theta_i  \cr}. $$
In  algebraic geometry, we
would regard $\theta_i$ and $ \pi_i$
as the homogeneous coordinates of a (complex) projective
line $\pp^1$, but in phylogenomics we limit
our attention to the real segment
specified by  $\theta_i \geq 0$, $\pi_i \geq 0 $ and
$\,\theta_i + 3  \pi_i = 1$.

Let $\Delta_{4^n-1}$ denote the set
of all probability distributions on $\Omega^n$.
Since $\Omega^n$ has $4^n$ elements, namely
the DNA sequences of length $n$, the set
$\Delta_{4^n-1}$ is a simplex of dimension
$4^n-1$. We identify the $j$-th leaf of
our tree $T$ with the \break $j$-th coordinate
of a DNA sequence $\,(u_1, \ldots, u_n) \in \Omega^n$,
and we introduce an unknown $\, p_{u_1 u_2 \cdots u_n}\,$  to represent
the probability of observing the nucleotides
$u_1,u_2,\ldots,u_n$ at the leaves $1,2,\ldots,n$. The
$4^n$ quantities $\, p_{u_1 u_2 \cdots u_n}\,$
are the coordinate functions on the simplex $\Delta_{4^n-1}$,
or, in the setting of algebraic geometry, on the projective
space $\,\pp^{4^n-1}\,$ obtained by complexifying $\Delta_{4^n-1}$.

\begin{prop} \label{JCprop}
In the Jukes-Cantor model on a tree $T$ with $r$ edges, the probability
     $\, p_{u_1 u_2 \cdots u_n}\,$ of making the
observation $(u_1,u_2,\ldots,u_n) \in \Omega^n$
at the leaves is expressed as a multilinear polynomial of degree $r$
in the model parameters
$\,(\theta_1,\pi_1), \,(\theta_2,\pi_2),
\ldots, \,(\theta_n,\pi_n) $.
Equivalently, in more geometric terms, the Jukes-Cantor 
model on $T$ is the image of a multilinear map
\begin{equation}
\label{JCmap}
     \phi \,\,\, : \,\,\,  (\pp^1)^r \,\,\longrightarrow \,\,  
\pp^{4^n-1}.
\end{equation}
\end{prop}

The coordinates of the map $\phi$ are easily
derived from the assumption
that the substitution processes along different edges
of $T$ are independent.
It turns out that the $4^n$ coordinates of $\phi$
are not all distinct.
To see this, we work out the formulas explicitly for a
very simple tree with three leaves.

\begin{ex} \label{jcdna}
Let $n=r=3$, and  let $T$ be the tree with
three leaves, labeled by $\, X = \{1,2,3\}$,
directly branching off the root of $T$.
We consider the Jukes-Cantor DNA model with
uniform root distribution on $T$. This model
is a three-dimensional algebraic variety, given
as the image of a trilinear map
$$ \phi \,\,:\,\, \pp^1 \times \pp^1 \times \pp^1 \,\,\,
\rightarrow \,\,\, \pp^{63}. $$
The number of states in $\Omega^3$ is $4^3 = 64$
but there are only five
distinct polynomials occurring among the
coordinates of the map $\phi$.
Let $p_{123}$ be the probability of observing the same letter at
all three leaves, $p_{ij}$ the probability of observing the  same
letter at the leaves $i,j$ and a different one at the third leaf,
and $p_{dis}$ the probability of seeing three distinct letters.
Then
\begin{eqnarray*}
     p_{123} \quad  = & \quad
\theta_1 \theta_2 \theta_3 \,\, + \,\,
3 \pi_1 \pi_2 \pi_3 , \\
     p_{dis} \quad   = & \quad
6 \theta_1 \pi_2 \pi_3 \,+\, 6 \pi_1 \theta_2 \pi_3 \,+\, 6 \pi_1
\pi_2 \theta_3 \,+\,
6 \pi_1 \pi_2 \pi_3 , \\
     p_{12} \quad  = & \quad
3 \theta_1 \theta_2 \pi_3 \,+ \, 3 \pi_1 \pi_2 \theta_3 \,+ \,
6 \pi_1 \pi_2 \pi_3 , \\
     p_{13} \quad  = & \quad
3 \theta_1 \pi_2 \theta_3  \,+ \, 3 \pi_1 \theta_2 \pi_3 \,+ \,
6 \pi_1 \pi_2 \pi_3 , \\
     p_{23} \quad  = & \quad
3 \pi_1 \theta_2 \theta_3 \,+ \, 3 \theta_1 \pi_2 \pi_3 \,+ \, 6
\pi_1 \pi_2 \pi_3 .
\end{eqnarray*}
All $64$ coordinates of $\phi$ are given by these five
trilinear polynomials, namely,
\begin{eqnarray*}
p_{AAA} \, = \,p_{CCC} \, = \, p_{GGG} \, = \, p_{TTT} \quad = \quad
\frac{1}{4} \cdot p_{123}  \\
p_{ACG} \, = \,p_{ACT} \, = \,\,\, \cdots \,\,\, = \, p_{GTC} \quad =
\quad
\frac{1}{24} \cdot p_{dis}  \\
p_{AAC} \, = \, p_{AAT} \, = \,\,\, \cdots \,\,\, = \, p_{TTG}
\quad = \quad \frac{1}{12} \cdot p_{12}  \\
p_{ACA} \, = \, p_{ATA} \, = \,\,\, \cdots \,\,\, = \, p_{TGT}
\quad = \quad \frac{1}{12} \cdot p_{13}  \\
p_{CAA} \, = \, p_{TAA} \, = \,\,\, \cdots \,\,\, = \, p_{GTT}
\quad = \quad \frac{1}{12} \cdot p_{23}
\end{eqnarray*}
This means that our Jukes-Cantor model is the image of the  simplified
map
$$ \phi' \, : \, \pp^1 \times \pp^1 \times \pp^1 \,\rightarrow \,
\pp^{4}, \,\, \bigl((\theta_1 ,  \pi_1),(\theta_2 ,  \pi_2),
(\theta_3 ,  \pi_3)\bigr)
\mapsto (p_{123},p_{dis},p_{12},p_{13},p_{23}). $$
In order to characterize the image of $\phi'$ algebraically,
we perform the following linear change of coordinates:
\begin{eqnarray*}
& q_{111} \, = \,
       p_{123} + \frac{1}{3} p_{dis} - \frac{1}{3} p_{12}
- \frac{1}{3} p_{13} - \frac{1}{3} p_{23}  \,= \, (\theta_1 -
\pi_1)(\theta_2 - \pi_2)(\theta_3 - \pi_3)
\\
& q_{110} \,\, = \,\,
     p_{123} - \frac{1}{3} p_{dis} + p_{12} - \frac{1}{3} p_{13}
- \frac{1}{3} p_{23}
      \,\, = \,\,
(\theta_1 - \pi_1)(\theta_2 - \pi_2)(\theta_3 + 3 \pi_3)
\\
& q_{101}\,\, = \,\,
     p_{123} - \frac{1}{3} p_{dis} - \frac{1}{3} p_{12}
+ p_{13} - \frac{1}{3} p_{23} \,\, = \,\, (\theta_1 -
\pi_1)(\theta_2 + 3 \pi_2)(\theta_3 - \pi_3)
\\
& q_{011} \,\, = \,\, p_{123}  - \frac{1}{3} p_{dis} - \frac{1}{3}
p_{12}
     - \frac{1}{3} p_{13} + p_{23}
      \,\, = \,\,
(\theta_1 + 3 \pi_1)(\theta_2 - \pi_2)(\theta_3 - \pi_3)
\\
& q_{000} \,\,\, = \,\,\, p_{123}  + p_{dis} +  p_{12} +  p_{13} +
p_{23}
      \,\,\, = \,\,\,
(\theta_1 + 3 \pi_1)(\theta_2 + 3 \pi_2)(\theta_3 + 3 \pi_3)
\end{eqnarray*}
This reveals that our model is the hypersurface in $\pp^4$
whose ideal equals
$$ I_T \quad = \quad \langle \, q_{000} q_{111}^2 \,
     - \, q_{011} q_{101} q_{110} \, \rangle $$
If we set $\,\theta_i = 1 - 3  \pi_i\,$
then we get the additional
constraint $\,q_{000}=1 $.
\qed
\end{ex}

The construction in this example generalizes to arbitrary
trees $T$.  There exists a change of coordinates, simultaneously
on the {\em parameter space} $\, (\pp^1)^r \,$
and on the {\em probability space} $ \, \pp^{4^n-1}$, such
that the map $\phi$ in (\ref{JCmap}) becomes a monomial map
in the new coordinates. This change of coordinates
is known as the {\em Fourier transform}
or as the {\em Hadamard conjugation}
(see \cite{ES, HenPen, SS, SSE}).

We regard the Jukes-Cantor DNA model on a tree $T$ with $n$ 
leaves and $r$ edges
as an algebraic variety of dimension $r$ in $\pp^{4^n-1}$, namely, it is
the image of the map (\ref{JCmap}).
Its homogeneous prime ideal $I_T$ is generated by differences of
monomials
$q^a - q^b$ in the Fourier coordinates. In the phylogenetics
literature (including the books \cite{Fel, SeSt}),  the polynomials in
the ideal $I_T$ are
known as {\em phylogenetic invariants} of the model.
The following result was shown in \cite{SS}.

\begin{thm} The ideal $I_T$ which defines
the Jukes-Cantor model on a binary tree $T$ is
generated by monomial differences
$q^a - q^b$ of degree at most three.
\end{thm}

It makes perfect sense to allow 
arbitrary distinct stochastic matrices
$P(t)$ on the edges of the tree $T$.
The resulting model is the {\em general Markov model}
on the tree $T$.  Allman and Rhodes \cite{AR1, AR2} determined
the complete system of phylogenetic  invariants for the general
Markov model
on a trivalent tree $T$.

An important problem in phylogenomics is to identify the maximum
likelihood branch lengths,
given a phylogenetic $X$-tree $T$, a rate matrix $Q$ and an
alignment of sequences. For the Jukes-Cantor DNA model
on three taxa, described in Example \ref{jcdna}, the exact ``analytic'' solution of this optimization
problem
leads to an algebraic equation of degree $23$.
See \cite[\S 6]{HKS} for details.

Let us instead consider the maximum likelihood
estimation problem in the
much simpler case of the Jukes-Cantor DNA model on two taxa.
Here the tree $T$ has only two leaves, labeled by $X = \{1,2\}$,
directly branching off the root of $T$. The model is given by
a surjective bilinear map
\begin{equation}
\label{mapfortwo} \phi \,\, : \,\, \pp^1 \times \pp^1 \,\,\rightarrow
\,\, \pp^1 \,
, \,\,\,\, ((\theta_1, \pi_1), (\theta_2 ,\pi_2))\,\,\mapsto
\,\,(\,p_{12} , p_{dis} \,).
\end{equation}
The coordinates of the map $\phi$ are
\begin{eqnarray*}
p_{12} \quad = & \quad \theta_1 \theta_2 \,\, + \,\, 3 \pi_1 \pi_2 ,\\
p_{dis} \quad = & \quad 3 \theta_1 \pi_2 \, + \, 3 \theta_2 \pi_1 \, +
\,6 \pi_1 \pi_2 .
\end{eqnarray*}
As before, we pass to affine coordinates by setting
$\,\theta_i = 1-3\pi_i\,$ for $i=1,2$.

One crucial difference between the model (\ref{mapfortwo})
and Example \ref{jcdna} is that the parameters in
    (\ref{mapfortwo}) are {\em not identifiable}.
Indeed, the inverse image of any point in $\pp^1$
under the map $\phi$ is a curve in $\pp^1 \! \times \! \pp^1$.
Suppose we are given data consisting of
two aligned DNA sequences of length $n$
where $k$ of the bases are different.
The corresponding point in $\pp^1$ is
$\, u = (n-k,k) $. The inverse image of $u$ under
the map $\phi$ is the curve in the
affine plane with the  equation
$$ 12 n \pi_1 \pi_2 \, - \, 3 n \pi_1 \,- \, 3 n \pi_2 \,+ \, k \quad =
\quad  0. $$
Every point $(\pi_1,\pi_2)$ on this curve is an
{\em exact fit} for the data $\, u = (n-k,k) $.
Hence this curve equals the set of all maximum likelihood
parameters for this model and the given data.
We rewrite the equation of the curve as follows:
\begin{equation}
\label{thefiber}
     (1 - 4 \pi_1)(1 - 4 \pi_2)  \quad = \quad   1 - \frac{4k}{3n}.
    \end{equation}
Recall from (\ref{branchlength}) that the branch length from the
root to leaf $i$ equals
$$ 3 \alpha_i t_i  \quad = \quad - \frac{1}{4} \cdot {\rm log} \,{\rm
det} \bigl( P_i(t) \bigr)
\quad = \quad - \frac{3}{4} \cdot {\rm log}(1 -  4 \pi_i) . $$
By taking logarithms on both sides of (\ref{thefiber}), we see that
the curve of all maximum likelihood parameters
becomes a line in the  branch length coordinates:
\begin{equation}
\label{thefiber2}
3 \alpha_1 t_1 \,\, + \,\,
3 \alpha_2 t_2 \quad = \quad
- \frac{3}{4} \cdot {\rm log} \bigl(    1 - \frac{4k}{3n} \bigr).
\end{equation}
The sum on the left hand side is the distance
from leaf $1$ to leaf $2$ in the tree $T$.
Our discussion of the two-taxa model leads
to the following
formula which known in evolutionary 
biology \cite{Fel} under the name {\em Jukes-Cantor correction}:

\begin{prop}
\label{JCcorrection}
Given an alignment
of two sequences of length $n$, with $k$ differences between the bases,
the ML estimate of the branch length equals 
\[ \delta_{12} \quad = \quad -\frac{3}{4} \cdot {\rm log} \left(1-\frac{4k}{3n} \right). \]
\end{prop}

There has been recent progress on solving the likelihood equations
exactly for small trees \cite{CHS, CKS, HKS,SY}. We believe that
these results will be useful in designing new algorithms
for computing maximum likelihood branch lengths, and to better
understand the mathematical properties of existing methods
(such as fastDNAml  \cite{OMHO})
    which are widely used by computational biologists.

It may also
be the case that $T$ is unknown, in which case the problem is not to
select a point on a variety, but to select from (exponentially many)
varieties. This  problem is discussed in the next
section.

The evolutionary models  discussed above do not allow for
insertion and deletion events.
They also assume  that sites evolve
independently. Although many widely
used models are based on these assumptions,
biological reality calls for models that include
insertion and deletion events \cite{HB},
site interactions \cite{SH},
and the flexibility
to allow for genome dynamics such as rearrangements.
Interested mathematicians will find a cornucopia of fascinating
research problems arising from such more refined 
evolutionary models.

\section{Phylogenetic Combinatorics}

Fix a set $X$ of $n$ taxa. A {\em dissimilarity map} on $X$ is a
function $\delta:X \times X \rightarrow {\mathbb R}$ such that
$\delta(x,x)=0$ and $\delta(x,y)=\delta(y,x)$.
The set of all dissimilarity maps on $X$ is a real vector space of
dimension
${n \choose 2}$ which we identify with $\rr^{n \choose 2}$.
A dissimilarity map $\,\delta \,$ is called a {\em metric on $X$} if
the triangle inequality holds:
$$ \delta(x,z) \,\,\leq \,\,\delta(x,y) + \delta(y,z) \qquad \hbox{for}
\,\,
x,y,z \in X. $$
The set of all metrics on $X$ is a full-dimensional convex polyhedral
cone in $\rr^{n \choose 2}$, called the {\em metric cone}.
Phylogenetic combinatorics is concerned with the study of
certain subsets of  the metric cone which are relevant for biology.
This field was pioneered in the 1980's by Andreas Dress
and his collaborators; see Dress' 1998 
ICM lecture  \cite{DreTer} and the references given there.

Let $T$ be a phylogenetic $X$-tree whose edges have specified lengths.
These lengths can be arbitrary non-negative real numbers.
The tree $T$ defines a metric $\delta_T$ on $X$
as follows: $\,\delta_T(x,y)  \,$ equals
the sum of the lengths of the edges on the unique path
in $T$ between the leaves labeled by $x$ and $y$.

The {\em space of $X$-trees} is the following subset of the metric cone:
\begin{equation} \label{TreeSpace}
\mathcal{T}_X \quad = \quad \bigl\{ \,\, \delta_T \,\,: \,\,
\hbox{$T$ is a phylogenetic $X$-tree} \,\bigr\} \quad \subset \quad
\rr^{n \choose 2}.
\end{equation}
Metric properties of the tree space $\mathcal{T}_X$
and its statistical and biological significance were  studied by
Billera, Holmes and Vogtmann  \cite{BHV}.
The following classical {\em Four Point
Condition} characterizes membership in the tree space:

\begin{thm} \label{4point}
A metric $\delta$ on $X$ lies in $\mathcal{T}_X$ if and only if,
for any four taxa $\,u,v,x,y \in X$, $\,\,\delta(u,v)+\delta(x,y)
\leq \max\{\delta(u,x)+\delta(v,y),\,\delta(u,y)+\delta(v,x)\}$.
\end{thm}

We refer to the book \cite{SeSt} for a proof of this
theorem and several variants.
To understand the structure of $\mathcal{T}_X$,
let us fix the combinatorial type of a trivalent tree $T$.
The number of choices of such trees is the {\em Schr\"oder number}
\begin{equation}
\label{schroder}
   (2n-5) !! \quad = \quad 1 \cdot 3 \cdot 5 \cdot \,\, \cdots \,\,
\cdot (2n-7) \cdot (2n - 5).
\end{equation}
Since $X$ has cardinality $n$,
the tree $T$ has $2n-3$ edges, and each of these edges
corresponds to a {\em split} $(A,B)$ of the set $X$ into two non-empty disjoint subsets $A $ and $B$.
Let ${\it Splits}(T)$ denote the collection of all $2n-3$
splits $(A,B)$ arising from $T$.

Each split $(A,B)$ defines a {\em split metric} $\,\delta_{(A,B)}\,$ on
$X$ as follows:
\begin{eqnarray*}
\delta_{(A,B)}(x,y) \,\,\, = \,\,\, 0 \quad &
\hbox{ if \ ($ x \in A$ and $y \in A$) \ or \
                    ($ x  \in B$ and $y \in B$), }\\
\delta_{(A,B)}(x,y) \,\,\, = \,\,\, 1 \quad &
\hbox{ if \ ($ x \in A$ and $y \in B$) \ or \
      ($ y \in A$ and $x \in B$). }
\end{eqnarray*}
The vectors $\,\bigl\{ \delta_{(A,B)} \,:\,
(A,B) \in {\it Splits}(T) \,\bigr\} \,$
are linearly independent in $\rr^{n \choose 2}$.
Their non-negative span is a cone $\,\mathcal{C}_T\,$ isomorphic to the
orthant $\, \rr_{\geq 0}^{2n-3}$.

\begin{prop}
The space  $\mathcal{T}_X$ of all $X$-trees is the union of the
$(2n-5)!!$  orthants $\,\mathcal{C}_T$. It is hence a
simplicial fan of pure dimension $2n-3$ in $\rr^{n \choose 2}$.
\end{prop}

The tree space $\mathcal{T}_X$ can be identified combinatorially
with  a simplicial complex of pure dimension
$2n-4$, to be denoted  $\widetilde{\mathcal{T}}_X$.
The vertices of $\widetilde{\mathcal{T}}_X$
are the $\,2^{n-1}-1\,$ splits of the
set $X$. We say that two splits $(A,B)$ and $(A',B')$ are
{\em compatible} if at least one of the four sets $ A \cap A'$,
$A \cap B'$, $B \cap A'$ and $B \cap B'$ is the empty set.
Here is a combinatorial characterization of the tree space:

\begin{prop}
\label{propcomp} A collection of splits of
the set $X$ forms a face in the
simplicial complex
    $\,\widetilde{\mathcal{T}}_X \,$ if and only if
    that collection is pairwise compatible.
\end{prop}

The {\em phylogenetics problem} is to reconstruct a tree $T$ from $n$
aligned sequences. In principle, one can select from evolutionary
models for all possible trees in order to find the maximum likelihood
fit. Even if the maximum likelihood problem can be solved for each
individual tree, this approach becomes infeasible  in practice
when $n$ increases, because
of the combinatorial explosion in the
number  (\ref{schroder}) of trees.
   A number of alternative approaches have been suggested
that attempt to find evolutionary models which fit {\em summaries} of
the data. They build on the
characterizations of trees given above.

{\em Distance-based methods} are based on the
observation that trees can be encoded by metrics satisfying the
Four Point Condition (Theorem \ref{4point}).
Starting from a multiple sequence alignment, one can produce a
dissimilarity map on the set $X$ of taxa by computing the maximum
likelihood
distance between every pair of taxa, using
Proposition \ref{JCcorrection}. The resulting dissimilarity map
$\delta$ is
typically not a tree metric, i.e., it does not
actually lie in the tree space $\mathcal{T}_X$.
What needs to be done  is to replace $\delta$
by a nearby tree metric $\delta_T \in \mathcal{T}_X$.

The method of choice for most
biologists is the {\em neighbor-joining algorithm},
which provides an easy-to-compute
map from the cone of all metrics onto
$\mathcal{T}_X$. The algorithm is based on the following
``cherry-picking theorem''
\cite{saitou,Studier}:

\begin{thm}
\label{simplecherry}
Let $\delta$ be a tree metric on $X$. For every pair $i,j \in X$ set
\begin{equation}
\label{njformula}
   Q_\delta(i,j)
\quad = \quad (n-2) \cdot \delta(i,j)
\, - \, \sum_{k \neq i} \delta(i,k)
\, - \, \sum_{k \neq j} \delta(j,k) .
\end{equation}
Then the pair $\,x,y \in X \,$ that minimizes
$Q_\delta(x,y)$ is a cherry in the tree, i.e.,
$x$ and $y$ are separated by only one internal vertex $z$ in the tree.
\end{thm}

Neighbor-joining works as follows. Starting from an arbitrary
metric $\delta$ on $n$ taxa, one sets up the $n \times n$-matrix
$\,Q_\delta\,$ whose $(i,j)$-entry is given by the formula
(\ref{njformula}), and one identifies the minimum off-diagonal entry
$Q_\delta(x,y)$. If $\delta$ were a tree metric then
the internal vertex $z$ which separates the leaves $x$ and $y$
would have the following distance from any other leaf $k$ in the tree:
\begin{equation}
\label{zkformula}
   \delta(z,k) \quad = \quad
\frac{1}{2}\bigl(\delta(x,k)+\delta(y,k)-\delta(x,y)\bigr).
\end{equation}
One now removes the taxa $x,y$ and replaces them by a new taxon
$z$ whose distance to the remaining $n-2$ taxa is given by
(\ref{zkformula}). This replaces the $n \times n$ matrix
$\,Q_\delta\,$ by an $(n-1) \times (n-1)$ matrix, and one
iterates the process.

This neighbor-joining algorithm recursively constructs a tree $T$
whose metric $\delta_T$ is reasonably close to the given metric
$\delta$. If $\delta$ is a tree metric then the method
is guaranteed to reconstruct the correct tree.
More generally, instead of estimating pairwise distances, one can
attempt to (more accurately) estimate the sum of the branch lengths of
subtrees of size $m \geq 3$.

We define an {\em $m$-dissimilarity map} on $X$ to be a
function $\,\delta:X^m  \rightarrow {\mathbb R}\,$
such that  $\, \delta(i_1,i_2,\ldots,i_m) \, = \,
   \delta(i_{\pi(1)},i_{\pi(2)},\ldots,i_{\pi(m)} )\, $
for all permutations $\pi$ on $\{1,\ldots,m\}$
and   $\, \delta(i_1,i_2,\ldots,i_m) \, = \,0 \,$
if  the taxa $i_1,i_2,\ldots,i_m$ are not distinct.
The set of all $m$-dissimilarity maps on $X$ is a real vector space of
dimension
${n \choose m}$ which we identify with $\rr^{n \choose m}$.
Every $X$-tree $T$ gives rise to an $m$-dissimilarity map
$\delta_T$ as follows. We define
$\delta_T(i_1,\ldots,i_m)$ to be the sum of all branch lengths in
the subtree of $T$ spanned by $\,i_1,\ldots,i_m \in X$.

The following theorem \cite{CL, PS3} is a generalization of
Theorem \ref{simplecherry}. It leads to a generalized neighbor-joining
algorithm which provides a better approximation of the maximum
likelihood tree and parameters:

\begin{thm}
\label{cherrypickingthm}
Let $T$ be an $X$-tree and $m <n =  |X|$.
For any $i,j \in X$ set
\[ Q_T(i,j) \,\, = \,\, \left( \frac{n-2}{m-1} \right) \sum_{ Y \in {X
\setminus \{i , j\} \choose m-2}} \!\!\!\! \delta_T(i,j,Y) \,\,
   -  \!\! \sum_{Y \in {X
\setminus \{i\} \choose m-1}} \!\!\! \delta_T (i,Y) \,\, - \!\!
\sum_{Y \in {X \setminus \{j\} \choose m-1}} \!\!\! \delta_T (j,Y). \]
Then the pair $\,x,y \in X \,$ that minimizes
$\,Q_T(x,y)\,$ is a cherry in the tree $T$.
\end{thm}

The subset of $\rr^{n \choose m}$ consisting of all
$m$-dissimilarity maps $\,\delta_T \,$ arising from
trees $T$ is a polyhedral space which is the image
of the tree space $\,\mathcal{T}_X \,$ under
a piecewise-linear map $\,\rr^{n \choose 2}
\rightarrow \rr^{n \choose m}$. We do not know
a simple characterization of this $m$-version of
tree-space which extends the Four Point Condition.

Here is another natural generalization of the space of
trees. Fix an $m$-dissimilarity map
$\,\delta:X^m  \rightarrow {\mathbb R}\,$ and
consider any $(m-2)$-element subset
$\, Y \in {X \choose m-2} $.
We get an induced dissimilarity map
$\,\delta/_{\! Y} \,$ on $\,X \backslash Y\,$ by setting
$$ \delta/_{\! Y} (i,j) \quad = \quad \delta(i,j,Y)
\quad \qquad \hbox{for all} \,\, i,j \in X \backslash Y . $$
We say that $\delta$ is an {\em $m$-tree} if
$\,\delta/_{\! Y} \,$ is a tree metric for
all  $\, Y \in {X \choose m-2} $. Thus,
by Theorem \ref{4point},
an $m$-dissimilarity map $\delta$  on $X$ is an $m$-tree if
$$\,\,\delta(i,j,Y)+\delta(k,l,Y)
\,\,\leq \,\,
\max\{\delta(i,k,Y)+\delta(j,l,Y),\,\delta(i,l,Y)+\delta(k,j,Y)\} $$
for all $\,Y \in {X \choose m-2} \,$ and
all $\,i,j,k,l \in X \backslash Y $.

Let  $\,Gr_{m,n}\,$ denote the subset of
$\, \rr^{n \choose m}\,$ consisting of all
$m$-trees. The space $Gr_{m,n}$ is a
polyhedral fan which is slightly larger than
the {\em tropical Grassmannian} studied in
\cite{SpStGrass}. 
For every $m$-tree$\,\delta \in Gr_{m,n}\,$
there is an $(m-1)$-dimensional tree-like space
whose ``leaves'' are the taxa in $X$.
This is the {\em tropical linear space} 
defined in \cite{Spe}.
This construction, which is described in
\cite[\S 6]{SpStGrass} and \cite[\S 3.5]{PSbook},
specializes to the construction of an $X$-tree $T$
from its metric $\delta_T$ when $m = 2$.
The study of $m$-trees and the tropical
Grassmannian was anticipated in
\cite{DreTer, DreWen}.  The Dress-Wenzel
theory of {\em matroids with coefficients} \cite{DreWen}
contains our $m$-trees as a special case. The space
$\, Gr_{m,n}\,$ of all $m$-trees  is discussed
in the context of buildings in \cite{DreTer}. Note that
the tree space $\mathcal{T}_X$
in (\ref{TreeSpace}) is precisely the
tropical Grassmannian $Gr_{2,n}$.

It is an open problem to find a natural
and easy-to-compute projection from
$\,\rr^{n \choose m}\,$ onto $\, Gr_{m,n}\,$
which generalizes the neighbor-joining method.
Such a variant of neighbor-joining would be likely
to have applications for more intricate  biological data
that are not easily explained by a tree model.
We close this section by discussing an example.

\begin{ex} \label{G36}
Fix a set of six taxa, $X = \{1,2,3,4,5,6\}$,
and let $m = 3$. The space of  $3$-dissimilarity
maps  on $X$ is identified with $\rr^{20}$.
An element $\delta \in \rr^{20}$ is a
$3$-tree if $\,\delta/_i \,$ is a tree metric
on $\, X \backslash \{i\}\,$ for all $i$.
Equivalently,
$$ \delta(i,j,k) +  \delta(i,l,m) \,\, \leq \,\,
\max \bigl\{ \delta(i,j,l) +  \delta(i,k,m) ,\,
\delta(i,j,m) +  \delta(i,k,l) \bigr\} $$
for all $i,j,k,l,m \in X$. The set $Gr_{3,6}$
of all $3$-trees is a $10$-dimensional
polyhedral fan. Each cone in this fan contains
the $6$-dimensional linear space $L$
consisting of all
$3$-dissimilarity maps of the particular form
$$ \delta(i,j,k) \, = \,  \omega_i +  \omega_j +
\omega_k \qquad \hbox{for some} \,\, \omega \in \rr^6. $$
The quotient $\,Gr_{3,6}/L \,$ is
a $4$-dimensional fan in the
$14$-dimensional real vectorspace $\rr^{20}/L$.
Let $\tilde Gr_{3,6}$ denote the three-dimensional
polyhedral complex obtained by intersecting
$\,Gr_{3,6}/L \,$ with a sphere around the origin
in $\rr^{20}/L$.

It was shown in \cite[\S 5]{SpStGrass} that
   $\,\tilde Gr_{3,6}\,$ is a  three-dimensional
simplicial complex consisting of
$65$ vertices, $550$ edges, $1395$ triangles
and $1035$ tetrahedra. Each of the
$1035$ tetrahedra
parameterizes six-tuples of tree metrics
$$\bigl(\, \delta/_{\! 1},\,
\delta/_{\! 2},\,
\delta/_{\! 3},\,
\delta/_{\! 4},\,
\delta/_{\! 5},\,
\delta/_{\! 6} \, \bigr), $$
where the tree topologies on five taxa are fixed.
The homology of the tropical Grassmannian   $\,\tilde Gr_{3,6}\,$
is concentrated in the top dimension
and is free abelian:
$$  H_3 \bigl(\tilde Gr_{3,6}, \zz \bigr) \quad =
\quad \zz^{126}. $$

If $T$ is an $X$-tree and $\delta_T$ the corresponding
$3$-dissimilarity map (as in Theorem \ref{cherrypickingthm})
then it is easy to check that $\delta_T$ lies in
$Gr_{3,6}$. The set of all $3$-trees
   $ \delta = \delta_T$ has codimension
one in $Gr_{3,6}$. It is the intersection
of  $Gr_{3,6}$ with the $15$-dimensional
linear subspace of $\rr^{20}$ defined by the equations
\begin{eqnarray*}
  \delta(1 2 3) + \delta(1 4 5) + \delta(2 4 6) + \delta(3 5 6)
&  = &
     \delta(1 2 4) + \delta(1 3 5) + \delta(2 3 6) + \delta(4 5 6) ,  \\
  \delta(1 2 3) + \delta(1 4 5) + \delta(3 4 6) + \delta(2 5 6)
&  = &
   \delta(1 3 4) + \delta(1 2 5) + \delta(2 3 6) + \delta(4 5 6) , \\
  \delta(1 2 3) + \delta(2 4 5) + \delta(1 4 6) + \delta(3 5 6)
& = &
    \delta(1 2 4) + \delta(2 3 5) + \delta(1 3 6) + \delta(4 5 6) ,  \\
  \delta(1 2 3) + \delta(3 4 5) + \delta(2 4 6) + \delta(1 5 6)
& = &
   \delta(2 3 4) + \delta(1 3 5) + \delta(1 2 6) + \delta(4 5 6) , \\
  \delta(1 2 3) + \delta(3 4 5) + \delta(1 4 6) + \delta(2 5 6)
&  = &
    \delta(1 3 4) + \delta(2 3 5) + \delta(1 2 6) + \delta(4 5 6) .
\end{eqnarray*}
Working modulo $L$ and intersecting with
a suitable sphere, the tree space $\widetilde{\mathcal{T}}_X$
is a two-dimensional simplicial complex,
consisting of $105 = 5 !!$ triangles.
To be precise, the simplicial complex in
Proposition \ref{propcomp} is the join of this
triangulated surface with the $5$-simplex on $X$.
Theorem  \ref{cherrypickingthm}
relates to the following
   geometric picture: the
triangulated surface
$\widetilde{\mathcal{T}}_X$
   sits inside the triangulated threefold
   $\tilde Gr_{3,6}$, namely,
as the solution set of the five equations. \qed
\end{ex}

\section{Back to the Data}

In Section 2, a conjecture was proposed based on our finding that
the ``meaning of life'' sequence (\ref{MOL}) is present (without  
mutations,
insertions or deletions) in  orthologous regions in ten vertebrate
genomes. In this section we explain how the various ideas
outlined throughout this paper can be used to estimate the probability
that such an extraordinary degree of conservation would occur by
chance. The mechanics of the calculation also provide a glimpse into
the types of processing and analyses that are performed in
computational biology. Two research papers
dealing with this subject matter are \cite{BPMSKMH, DEL}.

What we shall compute in this section is the
  probability under the Jukes-Cantor
model that a single ancestral base that is not under selection (and is
therefore free to mutate) is identical in the ten present day
vertebrates.

\smallskip

{\bf Step 1 (the genomes)}: The National Center for Biotechnology
Information (NCBI -- {\tt http://www.ncbi.nlm.nih.gov/}) maintains a public database called GENBANK which
contains all publicly available genome sequences from
around the world. Large sequencing centers that receive public funding
are generally required to deposit raw sequences
into this database within 24 hours of processing by sequencing
machines, and thus many automatic pipelines have been set up for generating and depositing sequences. The growth in GENBANK
has been spectacular. The database contained only $680,000$ base pairs
when it was started in 1982, and this number went up to $49$ million by
1990. There are currently $44$ billion base pairs of DNA in GENBANK.

The ten genomes of interest are not all complete, but are all
downloadable from GENBANK, either in pieces mapped to chromosomes
  (e.g.~for human) or  as collections of subsequences called {\em  
contigs}
(for less complete genomes).

\smallskip

{\bf Step 2 (annotation)}: In order to answer our question we need to
know where genes are in the genomes. Some genomes have annotations that
were derived experimentally, but {\em all} the genomes are annotated
using HMMs (Section 4) shortly after the release of the sequence. These
annotations are performed by centers such as at UC Santa Cruz ({\tt
http://genome.ucsc.edu/}) as well as by individual authors of programs.
It remains an open problem to accurately annotate
genomes. But HMM programs are quite good on average.
For example, typically $98$\% of coding bases are predicted correctly
to be in genes. On the other hand, boundaries of exons are often
misannotated: current state of the art methods only achieve
accuracies of about $80$\% \cite{AC}.

\smallskip

{\bf Step 3 (alignment)}: We start out by performing
  a genome alignment. Current methods for aligning whole genomes are all
based, to varying degrees, on the pair HMM ideas of Section 5. Although
in practice it is not possible to
align sequences containing billions or even millions of base pairs with hidden
Markov models, pair HMMs are subroutines of more complex  
alignment
strategies where smaller regions for alignment are initially identified
from the entire genomes by fast string matching algorithms \cite{BP}.
The ten vertebrate whole genome alignments which gave rise
to Conjecture \ref{MOLconj} are
accessible at \  {\tt http://bio.math.berkeley.edu/genomes/}.

\smallskip

{\bf Step 4 (finding neutral DNA)}: In order to compute the
probability that a certain subsequence is conserved between genomes,
it is necessary to estimate the {\em neutral rate of evolution}.
This is done by estimating parameters for an evolutionary model
of base pairs in the genome that are not under selection, and are  
therefore
free to mutate. Since neutral regions are difficult to identify a-priori, commonly used surrogates are synonymous substitutions in
codons (Section 3). Because synonymous substitutions do not change the
amino acids, it is unlikely that they are selected for or against, and
various studies have shown that such data provide good estimates
for neutral mutation rates. By searching through the annotations and
alignments, we identified $n = 14,202$ four-fold degenerate sites. These can be used for analyzing probabilities of neutral mutations.

\smallskip

{\bf Step 5 (deriving a metric)}:
We would ideally like to use maximum
likelihood techniques to reconstruct a tree $T$ with branch lengths
from the alignments of the four-fold degenerate sites.
One approach is to try to use a maximum-likelihood approach,
but this is difficult to do reliably because
of the complexity of the likelihood equations, even for the
Jukes-Cantor models with $|X| = 10$. An alternative approach is to
estimate
pairwise distances between species $i,j$ using the formula in
Proposition \ref{JCcorrection}. The resulting metric on the set
$\,X = \{{\rm gg}, {\rm hs}, {\rm mm}, {\rm pt} , {\rm rn}, {\rm  cf},
{\rm dr}, {\rm tn}, {\rm tr}, {\rm xt} \} \,$
is given in Table 3. For example, the pairwise alignment between
human and chicken (extracted from the multiple alignment) has
$n=14202$ positions, of which $k=7132$ are different.
Thus, the Jukes-Cantor distance
between the genomes of human and chicken equals
$$ -\frac{3}{4} \cdot {\rm log} \left(1-\frac{4k}{3n} \right)
\quad = \quad
  -\frac{3}{4} \cdot {\rm log} \left(\frac{14078}{42606}\right)
\quad = \quad 0.830536... $$

\begin{table}
\centering
\small
\begin{tabular}{|c|cccccccccc|}\hline
  & gg & hs & mm & pt & rn & cf & dr & tn & tr & xt\\
\hline
gg & -- & 0.831 &  0.928 & 0.831 & 0.925 & 0.847 &  1.321 & 1.326 &  
1.314 & 1.121\\
hs & -- & -- & 0.414 & 0.013 & 0.411 & 0.275 & 1.296 & 1.274 & 1.290 &  
1.166\\
mm & -- & -- & --    & 0.413 & 0.176 & 0.441 & 1.256 & 1.233 & 1.264 &  
1.218\\
pt & -- & -- & --    & --    & 0.411 & 0.275 & 1.291 & 1.267 & 1.288 &  
1.160\\
rn & -- & -- & --    & --    &--     & 0.443 & 1.255 & 1.233 & 1.258 &  
1.212\\
cf & -- & -- & --    & --    & --    & --    & 1.300 & 1.251 & 1.269 &  
1.154\\
dr & -- & -- & --    & --    & --    & --    & --    & 1.056 & 1.067 &  
1.348\\
tn & -- & -- & --    & --    & --    & --    & --    & --    & 0.315 &  
1.456\\
tr & -- & -- & --    & --    & --    & --    & --    & --    & --    &  
1.437\\
\hline
\end{tabular}
\caption{Jukes-Cantor pairwise distance estimates.}
\end{table}

\begin{figure}[ht]
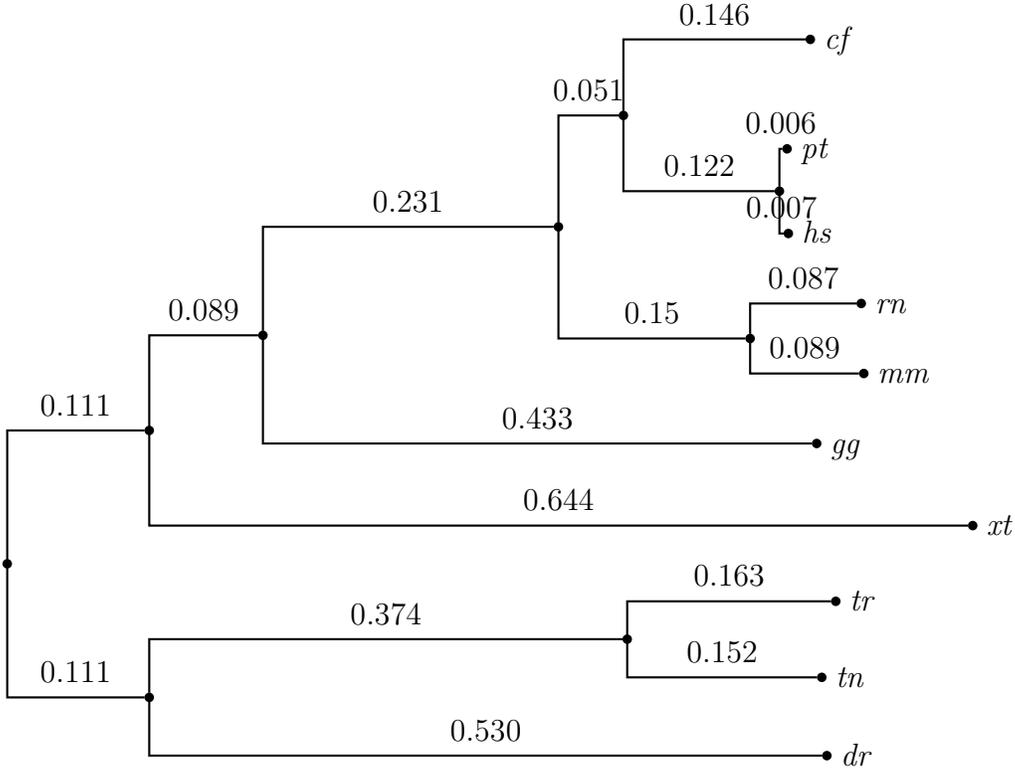

 \begin{newicktree}
 \setunitlength{17cm}
 \righttree
\drawtree{(((((cf:0.146,(pt:0.006,hs:0.007):0.122):0.051,(rn:0.087,mm:0.089):0.15):0.231,gg:0.433):0.089,xt:0.644):0.111,((tr:0.163,tn:0.152):0.374,dr:0.530):0.111);}
\end{newicktree}
\caption{Neighbor joining tree from alignment of codons in ten vertebrates.}
\label{NJtree}
\end{figure}

\smallskip
{\bf Step 6 (building a tree)}:
{}From the pairwise distances in Table 3 we construct a
phylogenetic $X$-tree using the neighbor joining algorithm (Section 7).
The tree with the inferred branch lengths is shown in Figure
\ref{NJtree}. The tree is drawn such that the branch lengths
are consistent with the horizontal distances in the diagram.
 The root of the tree was added manually in order to properly
indicate the ancestral relationships between the species. 

At this point we wish to add a philosophical remark:
{\sl The tree in Figure \ref{NJtree} is a point on an algebraic
variety!}
Indeed, that variety is the Jukes-Cantor model (Proposition
\ref{JCprop}),
and the preimage coordinates $(\theta_i, \pi_i)$
of that point are obtained by exponentiating the branch lengths
as described in Section 6.

\smallskip

{\bf Step 7 (calculating the probability)}:
We are now given a specific point on the variety
representing the Jukes-Cantor model on the tree
depicted in Figure \ref{NJtree}. Recall from
Proposition   \ref{JCprop} that this variety,
and hence our point, lives in a projective space of dimension
$\, 4^{10}-1 =  1,048,575 $.  What we are interested in
are four specific coordinates of that point, namely,
the probabilities
that the same nucleotide occurs in every species:
\begin{equation}
\label{AllSame}
p_{AAAAAAAAAA}\, = \,
p_{CCCCCCCCCC} \,= \,
p_{GGGGGGGGGG} \,= \,
p_{TTTTTTTTTT}
\end{equation}

As discussed in Section 6, this expression
is a multilinear polynomial in the edge
parameters $(\theta_i,\phi_i)$. When we
evaluate it at the parameters derived from
the branch lengths in Figure \ref{NJtree}
we find that
$$ p_{AAAAAAAAAA} \quad  = \quad  0.009651... $$
Returning to  the ``meaning of life'' sequence (\ref{MOL}),
this implies the following

\begin{prop} \label{WrapItUp}
Assuming the probability distribution on $\Omega^{10}$ given
by the Jukes-Cantor model on the tree in Figure \ref{NJtree},
the probability of observing a sequence of length $42$
unchanged at a given location
in the ten vertebrate genomes within a neutrally
evolving region equals $ \,  (0.038604)^{42} \,\,  = \,\, 4.3 \cdot 10^{-60} $.
\end{prop}

This calculation did not take into account the fact that
the ``meaning of life'' sequence may occur in an arbitrary location of
the genome in question. In order to adjust for this, we can multiply the number
 in Proposition \ref{WrapItUp} by the length
of the genomes. The human genome contains approximately
$2.8$ billion nucleotides, so it is reasonable to conclude that
the probability of observing a sequence of length $42$
unchanged {\em somewhere}  in the ten vertebrate genomes
is approximately
$$ 2.8 \cdot 10^{9} \,\, \times \,\, 4.3 \cdot 10^{-60} \quad \simeq \quad 10^{-50} .$$
This probability is a very small number, i.e.,  it is unlikely
that the remarkable properties of the
 sequence (\ref{MOL}) occurred by ``chance''.
Despite the shortcomings of the Jukes-Cantor model
discussed at the end of Section 6, we believe that
Proposition \ref{WrapItUp} constitutes a sound argument in support of Conjecture~1.

\section{Acknowledgments}

The vertebrate whole genome alignments we have
analyzed were assembled by Nicolas Bray and Colin Dewey. We also thank  
Sourav
Chatterji and Von Bing Yap for their help in searching through the  
alignments.
Lior Pachter was supported by a grant from the NIH (R01-HG2362-3), a Sloan Foundation Research Fellowship, and an NSF CAREER award  (CCF-0347992). Bernd Sturmfels was supported by the NSF (DMS-0200729, DMS-0456960).


\begin{thebibliography}{99}


\bibitem{A} D.~N.~Adams. {\em The Hitchhikers Guide to the Galaxy}, Pan
Books, 1979.

\bibitem{ABP} M.~Alexandersson, N.~Bray and L.~Pachter. Pair hidden Markov models. Encyclopedia of Genetics, Genomics, Proteomics and Bioinformatics (L.~B.~Jorde, P.~Little, M.~Dunn and S.~Subramanian, editors), in press.

\bibitem{ACP} M.~Alexandersson, S.~Cawley and L.~Pachter. SLAM--
Cross-species gene finding and alignment with a generalized pair hidden
Markov model. {\em Genome Research} 13 (2003) 496--502.

\bibitem{AR1} E.~Allman and J.~Rhodes.  Phylogenetic  invariants for
     the general Markov model of sequence  mutation.
     \emph{Mathem.~Biosciences} 186 (2003) 133--144.

\bibitem{AR2} E.~Allman and J.~Rhodes.      
Phylogenetic ideals and varieties for the general Markov model,
preprint, {\tt math.AG/0410604}.

\bibitem{AC} J.~Ashurst and J.~E.~Collins. Gene annotation: prediction
and testing, {\em Annual Review of Genomics and Human Genetics} 4
(2003) 69--88.

\bibitem{BPMSKMH} G.~Bejerano, M.~Pheasant, I.~Makunin, S.~Stephen,
W.~J.~Kent, J.~S.~Mattick and D.~Haussler. Ultraconserved elements in
the human genome, {\em Science} 304 (2004) 1321--1325.

\bibitem{Bic} P.~Bickel and K.~Doksum.
{\em Mathematical Statistics: Basic Ideas and Selected Topics},  
Holden-Day Inc., San Francisco, Calif., 1976.

\bibitem{BHV} L.~Billera, S.~Holmes, and K.~Vogtmann.
Geometry of the space of phylogenetic trees,
{\em Advances in Applied Mathematics} 27 (2001) 733--767.

\bibitem{BP} N.~Bray and L.~Pachter. MAVID: Constrained ancestral
alignment of multiple sequences, {\em Genome Research} 14 (2004)
693--699.

\bibitem{BH} P.~Bucher and K.~Hofmann. A sequence similarity search
algorithm based on a probabilistic interpretation of an alignment
scoring system, {\em Proceedings of ISMB} (1996) 44--51.

\bibitem{BK} C.~Burge and S.~Karlin. Prediction of complete gene
structures in human genomic DNA, {\em Journal of Molecular Biology} 268
(1997) 78--94.

\bibitem{CMK} A.~Campbell, J.~Mrazek and S.~Karlin. Genome signature
comparisons among prokaryote, plasmid and mitochondrial DNA.
 {\em Proc. Natl. Acad. Sci. USA} 96 (1999) 9184--9189.
\bibitem{CP} S.~Chatterji and L.~Pachter. Multiple organism gene
finding by collapsed Gibbs sampling, {\em Proceedings of the Eighth
Annual International Conference on Computational Molecular Biology
-- RECOMB 2004}, San Diego, April 2004, pp.~187--193.

\bibitem{CHS} B.~Chor, M.~Hendy and S.~Snir. Maximum likelihood
Jukes-Cantor triplets: analytic solutions,
preprint, {\tt ArXiv:q-bio.PE/0505054}.

\bibitem{CKS} B.~Chor, A.~Khetan and S.~Snir.
Maximum likelihood on four taxa phylogenetic trees: analytic solutions,
{\em Proceedings of the Seventh Annual Conference on Research in  
Computational
Molecular Biology -- RECOMB 2003}, Berlin, April 2003, pp.~76--83.

\bibitem{CL} M.~Contois and D.~Levy. Small trees and generalized neighbor-joining.
In \cite{PSbook}, pp.~335--346.

\bibitem{DEKM} R.~Durbin, S.~R.~Eddy, A.~Korgh and G.~Mitchison. Biological Sequence Analysis: Probabilistic
Models of Proteins and Nucleic Acids, Cambridge University Press, 1999.

\bibitem{DreTer}
A.~Dress and W.~Terhalle.
The tree of life and other affine buildings.
Proceedings of the International Congress of Mathematicians,
Vol. III (Berlin, 1998), {\em Documenta Mathematica}
(1998) Extra Vol. III, 565--574.

\bibitem{DreWen}
A.~Dress and W.~Wenzel.
Grassmann-Pl\"ucker relations and matroids with coefficients,
{\em Advances in  Mathematics} 86 (1991)   68--110.

\bibitem{DEL}
M.~Drton, N.~Eriksson and G.~Leung.
Ultra-conserved elements in vertebrate and fly genomes,
in \cite{PSbook},  pp.~387--402.

\bibitem{E}
J.~Eisen.
Phylogenomics: Improving functional predictions for uncharacterized genes by evolutionary analysis, {\it Genome Research} 8 (1998) 163--167.

\bibitem{EichS} E.~E.~Eichler and D.~Sankoff. Structural dynamics of  
eukaryotic chromosome evolution, {\it Science} 301 (2003) 793--797.

\bibitem{Eli} S.~Elizalde. Inference functions.
In \cite{PSbook}, pp.~215--225.


\bibitem{ES} S.~Evans and T.~Speed.  Invariants of some probability models used in phylogenetic inference, {\em Annals of Statistics}
21 (1993) 355--377.

\bibitem{Fel} J.~Felsenstein.  \emph{Inferring Phylogenies}.  Sinauer
     Associates, Inc., Sunderland, 2003.

\bibitem{FV} D.~Fern'{a}ndez-Baca and B.~Venkatachalam. Parametric
Sequence Alignment, {\em Handbook on Computational Molecular Biology}
(S. Aluru, ed.), Chapman and Hall/CRC press,  2005.

\bibitem{Fie} S.E.~Fienberg. {\em The Analysis of Cross-classified Categorical
Data},  2nd edition, M.I.T. Press, Cambridge, MA, 1980.

\bibitem{Gusfield:97} D.~Gusfield. {\em Algorithms on Strings, Trees and
Sequences}, Cambridge University Press, 1997.

\bibitem{HenPen}
M.D.~Hendy and D.~Penny. Spectral analysis of phylogenetic data,
{\em Journal of Classification} 10 (1993) 5--24.

\bibitem{HB} I.~Holmes and W.~J.~Bruno. Evolutionary HMMs: A bayesian
approach to multiple alignment,  {\em Bioinformatics} 17(9) (2001) 803--820.

\bibitem{HKS} S.~Ho\c{s}ten, A.~Khetan and B.~Sturmfels. Solving the
likelihood equations, {\em Foundational of Computational
Mathematics}, to appear. 

\bibitem{HMY} I.~Hallgr\'{\i}msd\'{o}ttir, R.~A.~Milowski and J.~Yu. The EM algorithm for hidden Markov models. In
\cite{PSbook}, pp 250--263.

\bibitem{DiscreteThoughts} M.~Kac, G-C.~Rota and J.~T.~Schwartz. {\em Discrete Thoughts}, Birkh\"{a}user, Boston,1986.

\bibitem{Karp} R.M.~Karp. Mathematical challenges from genomics and
molecular biology, {\em Notices of the American Mathematical Society}
49 (2002) 544--553.

\bibitem{Kellis} M.~Kellis, B.~Birren and E.~Lander. Proof and
evolutionary analysis of ancient genome duplication in the yeast
Saccharomyces cerevisiae, {\em Nature} 8 (2004) 617--624.

\bibitem{KHRE} D.~Kulp, D.~Haussler, M.G.~Reese, and F.H.~Eeckman. A
generalized hidden Markov model for the recognition of human genes in
DNA, {\em Proceedings of ISMB} 4 (1996) 134--142.

\bibitem{Lander} E.~S.~Lander et al. Initial sequencing and analysis of  
the human genome, {\it Nature} 409 (2001) 860--921.


\bibitem{Myers} E.~Myers et al. A whole-genome assembly of Drosophila,  
{\em Science} 287 (2000) 2196--2204.

\bibitem{OMHO} G.~J.~Olsen, H.~Matsuda, R.~Hagstrom and R.~Overbeek.
fastDNAml: A tool for construction of phylogenetic trees of DNA
sequences using maximum likelihood, {\em Comput. Appl. Biosci.} 10
(1994) 41--48

\bibitem{PS3} L.~Pachter and D.~Speyer.
Reconstructing trees from subtree weights,
{\em Applied Mathematics Letters} 17 (2004) 615--621.

\bibitem{PS1} L.~Pachter and B.~Sturmfels.
Tropical geometry of statistical models, 
  {\em Proc. Natl. Acad. Sci. USA} 101 (2004) 16132-16137

\bibitem{PS2} L.~Pachter and B.~Sturmfels.
    Parametric inference for biological sequence analysis,
  {\em Proc. Natl. Acad. Sci. USA} 10 (2004) 16138-16143. 

\bibitem{PSbook}  L.~Pachter and B.~Sturmfels.
{\em Algebraic Statistics for Computational Biology},
Cambridge University Press, 2005.

\bibitem{PT} P.~Pevzner and G.~Tesler. Human and mouse genomic  
sequences reveal extensive breakpoint reuse in mammalian evolution,
  {\em Proc. Natl. Acad. Sci. USA}  100 (2003) 7672--7677.

\bibitem{saitou} N.~Saitou and M.~Nei. The neighbor joining method: a
new method for reconstructing phylogenetic trees, {\em Molecular
Biology and Evolution} 4 (1987) 406--425.

\bibitem{SY} R. Sainudiin and R.~Yoshida. Applications
of interval methods to phylogenetics,
in \cite{PSbook}, pp.~359--374.

\bibitem{SN} D.~Sankoff and J.~H.~Nadeau. Chromosome rearrangements
in evolution: From gene order to genome sequence and back,
   {\em Proc. Natl. Acad. Sci. USA} 100 (2003) 11188--11189.

\bibitem{SeSt} C.~Semple and M.~Steel.  \emph{Phylogenetics}.
Oxford University Press, 2003.

\bibitem{SH} A. Siepel and D. Haussler. Phylogenetic estimation 
of  context-dependent substitution rates by maximum likelihood.  
{\em Molecular Biology and Evolution} 21 (2004) 468-488.

\bibitem{Spe} D.~Speyer. Tropical linear spaces,
preprint, {\tt math.CO/0410455}.

\bibitem{SpStGrass} D.~Speyer and B.~Sturmfels.
The tropical Grassmannian, {\sl Advances in Geometry}
{\bf 4} (2004) 389--411.

\bibitem{Stanley:99} R.P.~Stanley. {\em Enumerative Combinatorics},  
Vol. 1,
Cambridge Studies in Advanced Mathematics, {\bf 49}, Cambridge
University Press, 1997.

\bibitem{Studier} J.~A.~Studier and K.~J.~Keppler. A note on the
neighbor-joining method of Saitou and Nei, {\em Molecular
Biology and Evolution} 5(1988) 729--731.

\bibitem{SS} B.~Sturmfels and S.~Sullivant.
Toric ideals of phylogenetic invariants, {\em Journal of Computational  
Biology} 12 (2005) 204-228. 

\bibitem{SSE} L.~Sz\'ekely, M.~Steel and P.~Erd\"os. Fourier calculus
on evolutionary trees, {\em Advances in Applied Mathematics}
14 (1993) 200-210.

\bibitem{Venter} J.~C.~Venter et al. The sequence of the human genome,
{\em Science} 291 (2001) 1304--1351.

\bibitem{Waterman:92} M.~Waterman, M.~Eggert and E.~Lander. Parametric
sequence comparisons, {\em Proc. Natl. Acad. Sci. USA} 89 (1992)
6090--6093.

\bibitem{WC} J.~Watson and F.~Crick. A structure for Deoxyribose  
Nucleic Acid, {\em Nature} 171 (1953) 964-967.

\bibitem{YP} V.B.~Yap and L.~Pachter. Identification of evolutionary
hotspots in the rodent genomes, {\em Genome Research} 14 (2004)  
574--579.
    \end{thebibliography}
\end{document}